\documentclass[11pt]{amsart}
\usepackage{amssymb,amsfonts}
\usepackage{amsmath,amsthm}
\usepackage{color}
\newcommand{\bg}{\begin{equation}}
\newcommand{\ed}{\end{equation}}
\newcommand{\bga}{\begin{eqnarray}}
\newcommand{\eda}{\end{eqnarray}}

\def\ep{\varepsilon}
\def\epunu{{\varepsilon_1}}
\def\epdoi{{\varepsilon_2}}
\def\eptrei{{\varepsilon_3}}
\def\eppatru{{\varepsilon_4}}
\def\ub{{\overline u}}
\def\vb{{\overline v}}
\def\qb{{q_1}}
\DeclareMathOperator{\mes}{mes}

\newtheorem{theorem}{Theorem}[section]
\newtheorem{proposition}[theorem]{Proposition}
\newtheorem{lemma}[theorem]{Lemma}

\theoremstyle{definition}

\theoremstyle{remark}
\newtheorem{remark}[theorem]{Remark}

\numberwithin{equation}{section}

\def\R{\mathbb{R}}
\def\C{\mathbb{C}}

\def\endProof{{\hfill$\Box$}}

\renewcommand{\P}{\mathbb{P}}

\begin{document}

\title[$L^p$-solutions of the steady-state Navier--Stokes equations (\today)]
{
$L^p$-solutions of the steady-state Navier--Stokes\\
equations with rough external forces}

\author[Bjorland]{Clayton Bjorland}
\author[Brandolese]{Lorenzo Brandolese}

\author[Iftimie]{Drago\c s Iftimie}

\author[Schonbek]{Maria E. Schonbek}
\address{C. Bjorland: Department of Mathematics, UC Santa Cruz, Santa Cruz, CA 95064,USA}
\email{cbjorland@math.ucsc.edu}
\address{L. Brandolese: Universit\'e de Lyon~; Universit\'e Lyon 1~;
CNRS UMR 5208 Institut Camille Jordan,
43 bd. du 11 novembre,
Villeurbanne Cedex F-69622, France.}
\email{brandolese{@}math.univ-lyon1.fr}
\urladdr{http://math.univ-lyon1.fr/\~{}brandolese}

\address{D. Iftimie: Universit\'e de Lyon~; Universit\'e Lyon 1~;
CNRS UMR 5208 Institut Camille Jordan,
43 bd. du 11 novembre,
Villeurbanne Cedex F-69622, France.}
\email{iftimie{@}math.univ-lyon1.fr}
\urladdr{http://math.univ-lyon1.fr/\~{}iftimie}

\address{M Schonbek: Department of Mathematics, UC Santa Cruz, Santa Cruz, CA 95064,USA}
\email{schonbek@math.ucsc.edu}

\thanks{ The work of L. Brandolese, C. Bjorland, D. Iftimie and M. Schonbek were partially supported by FBF GrantSC-08-34.
The work of M. Schonbek was also partially supported by NSF Grant DMS-0600692.}


\date{\today}


\begin{abstract}
 In this paper we address the existence, the asymptotic behavior
 and stability in $L^p$ and $L^{p,\infty}$, $\frac{3}{2}<p\le\infty$,
for solutions to the steady state 3D Navier-Stokes equations with possibly very singular
external forces.
We show that under certain smallness conditions of the forcing term
 there exists solutions to the stationary Navier-Stokes equations in $L^p$ spaces,
 and we prove the stability of these solutions.
Namely, we prove that such small steady state solutions attract time dependent solutions driven by the same
forcing, no matter how large the initial velocity of non-stationary solutions is.

We also give non-existence results of  stationary solutions in $L^p$, for $1\le p\le \frac{3}{2}$.
\end{abstract}

\maketitle

\section{Introduction}
In this paper we consider the solutions to the three-dimensional steady state Navier--Stokes equations in the whole space
$\R^3$,

\begin{equation}
\label{SNS} 
\left\{
\begin{aligned}
 &\nabla\cdot (U\otimes U)+\nabla P=\Delta U+f\\
 &\nabla\cdot U=0.\\
\end{aligned}
\right.
\end{equation}
Here $U=(U_1,U_2,U_3)$ is the velocity, $P$ the pressure and $f=(f_1,f_2,f_3)$
a given time independent external force.
Equation  \eqref{SNS}  will be complemented with a boundary condition at infinity
of the form~$U(x)\to 0$  in a weak sense: typically, we express
this condition requiring that~$U$ belongs to some $L^p$ spaces.
Three problems will be addressed. 

We will first establish the existence of solutions  $U\in L^p$,
with $\frac{3}{2}<p\le \infty$,  to equations (\ref{SNS})
for (small) functions $f$  as general as possible, and non-existence results in the range $1\le p\le \frac{3}{2}$.

Next we will study the asymptotic properties as $|x|\to\infty$
for a relevant subclass of the solutions obtained.

The third problem at hand is the stability of the solutions in the sense of  solutions to (\ref{SNS}) being ``fixed point''  in $L^p$
to the non-stationary  incompressible Navier-Stokes equations in $\R^3$

\begin{equation}\label{NS:PDE}
\left\{
\begin{aligned}
&\partial_tu+u\cdot\nabla u+\nabla p = \Delta u +f\\
&\nabla \cdot u =0\\
&u(0)=u_0,
\end{aligned}
\right.
\end{equation}

\noindent
where  $u$, $p$ are the time dependent velocity and pressure of the flow.
We assume $f$  to be constant in time, but our methods could also be applied
to the more general case of time dependent forces suitably converging to a steady state forcing term.
We will show that  \emph{small} stationary solutions $U$ of~\eqref{SNS} will attract
all global non-stationary solutions~$u$ to~\eqref{NS:PDE} verifying mild regularity conditions,
and emanating from possibly \emph{large} data~$u_0$.
This will be achieved by first proving that a wide class of global solutions  of \eqref{NS:PDE} must become small in $L^{3,\infty}$
after some time, and then applying the  stability theory of small solutions in~$L^{3,\infty}$ as developed, e.g.,
in~\cite{CanK05, KozY98, Yam00}.
In addition, for small solutions, we will extend the results on the stability in the existing literature by giving
necessary and sufficient conditions to have $u(t)\to U$ in $L^p$ as $t\to\infty$.

\medskip
 The existence and stability of stationary solutions  is well understood in the case of bounded domains.
 See for example \cite{DG}.
 For related results  in exterior domains we refer the reader to
 \cite{Finn59,Finn61,Finn65,H}. 
A wider list of references regarding connected literature  can be found in \cite{BjorSch07}.
For example,
the existence and the stability  of stationary solutions in  $L^p$  with $p \geq n$,
where $n$ is the dimension of the space,
is obtained in \cite{Secchi}, under the condition that
the  Reynolds number is sufficiently small,
and in \cite{KozY98}, \cite{Yam00} under the assumption that
the external force is small in a Lorentz space.
Similar results in the whole domain $\R^n$, always for $p\ge n$,
have been obtained also in~\cite{KozY95}, \cite{CanK}, \cite{CanK05}.

On the other hand, not so much can be found in the literature about
the existence and stability of stationary solutions in $\R^n$
with $p<n$.
This problem have been studied recently in the case $n=3$ and $p=2$ 
in \cite{BjorSch07}.
In this paper
we extend the results of \cite{BjorSch07} to the range $\frac{3}{2}<p\le \infty$,
and improve such results also in the case $p=2$ by considering a more general
class of forcing functions. The methods in this paper   differ completely from the ones used in 
\cite{BjorSch07}. In the former paper the construction of solutions with finite energy was based on a well known formal observation: if $\Phi$ is the fundamental solution for the heat equation then $\int_0^\infty \Phi(t,\cdot)\,dt$ is the fundamental solution for Poisson's equation.  Using that idea it was possible to make a time dependent PDE similar to the Navier-Stokes equation with $f$ as initial data with a solution that can be formally integrated in time to find a solution of (\ref{SNS}). 

As we shall see, the conditions on~$f$ in the present paper which  yield that  $U\in L^p$ are, essentially,
necessary and sufficient.
This will be made possible by a systematic use of suitable function spaces.

\medskip
One could also complement the system~\eqref{SNS} with different type of boundary condition at infinity.
For example, conditions of the form
 $U(x)\to U_\infty$ as $|x|\to\infty$,
where $U_\infty\in\R^3$ and $U_\infty\not=0$ are also of interest.
However the properties of stationary solutions satisfying such condition are  already quite well understood.
We refer to the treatise of Galdi~\cite{GaldiBook}
for a comprehensive study of this question.

On the other hand, the understanding of the problem in the case $U_\infty=0$
is less satisfactory.
For example, the construction of solutions obeying to the natural energy equality
(obtained multiplying the equation~\eqref{SNS} by $U$ and formally integrating by parts),
without putting any smallness assumption on~$f$, is still an open problem.
The main difficulty, for example when $\Omega=\R^3$ (or when  Poincar\'e's inequality is not available),
is that the usual {\it a priori\/} estimate on the Dirichlet integral
$$\|\nabla U\|_{L^2}\le \|f\|_{\dot H^{-1}}$$
ensures only that $U\in \dot H^1\subset L^6$:
but to give a sense to the integral in the formal equality 
$$\int \bigl[\nabla\cdot (U\otimes U)\bigr]\cdot U\,dx=0$$
one would need, {\it e.g.\/}, that $U$ belongs also to $L^4$.

More generally, one motivation for developing the $L^p$ theory
 (especially for {\it low values\/} of~$p$) of stationary solutions  is that 
this provides additional information on the asymptotic properties of~$U$ in the far field.
On the other hand,  condition like $U\in L^p$ for large~$p$ are usually easily recovered
via the standard regularity theory, as bootstrapping procedures show that
weak solutions $U\in \dot H^{1}$ are regular if $f$ is so. See also \cite{Secchi} for this case.

\medskip
The paper will be organized as follows. After the introduction we have a section of general notation, where we recall definitions of  several function
spaces which will be needed in the sequel.

Section two deals with the existence of solutions in $L^p$, $\frac{3}{2}<p\le \infty$.
Section three  addresses the pointwise behavior in $\R^3$ of the solutions and the asymptotic profiles. We note that the study of the
asymptotic profiles has been largely dealt  in the literature, starting with the well known results of Finn \cite{Finn65}
in exterior domains. Our results being in the whole domain are simpler, but we are able to get them  with weaker conditions.
Non existence results of (generic) solutions in $U\in L^p$, $p\le \frac{3}{2}$ will also follow from such analysis.

Section four  handles the stability of stationary solutions. More precisely in the setting of the 
Navier--Stokes equation we investigate the stability  of the stationary solution  $U$  in the $L^p$ and the Lorentz $L^{p,\infty}$-norms. We consider a possibly large $L^{3,\infty}$ non-stationary solution and a stationary solution $U\in L^{3,\infty}\cap L^p$ or $U\in L^{3,\infty}\cap L^{p,\infty}$ which is small in $L^{3,\infty}$. We show that the non-stationary solution eventually becomes small in $L^{3,\infty}$ (but does not converge to $0$ in this space), we prove some decay estimates for it and we give a necessary and sufficient condition to have that $u(t)\to U$ in $L^p$ or $L^{p,\infty}$.

The fact that small steady state solutions~$U$ attract small non-stationary solutions was proved
by several authors in different functional settings, see, {\it e.g.\/} \cite{CanK, CanK05, KozY95, KozY98, Yam00}. 
The main novelty of our approach is that we can prove the same result for a class of large solutions.
At best of our knowledge, this was known only in the particular case $U=0$ (see \cite{ADT, GIP}).
Our main tool will be a decomposition criterion for functions in Lorentz-spaces.

\subsection{Notations}

\subsubsection{Function spaces}

We recall that the fractional Sobolev spaces (or Bessel potential spaces) are
defined, for $s\in\R$ and $1<p<\infty$, as
\begin{equation*}
H^s_p=\{ f\in \mathcal{S}'(\R^3)\colon \mathcal{F}^{-1}(1+|\xi|^2)^{\frac s2}\widehat f\in L^p\},
\end{equation*}
and their homogeneous counterpart is
\begin{equation*}
\dot H^s_p=\{ f\in \mathcal{S}'(\R^3)\colon \mathcal{F}^{-1}|\xi|^{s}\widehat f\in L^p\}.
\end{equation*}
Their differential dimension is~$s-\frac3p$. We will only deal with the case $s-\frac3p<0$, so that
the elements of $\dot H^s_p$ can indeed be realized as tempered distributions.
As usual, we will simply write $H^s$ and $\dot H^s$ instead of $H^s_2$ and $\dot H^s_2$
for the classical Sobolev spaces.

The fractional Sobolev spaces can be identified with 
particular Triebel-Lizorkin spaces, namely $F^{s,2}_p$ and $\dot F^{s,2}_p$.
This identification will be useful, because it allows us to handle the limit case for $p=1$:
the corresponding spaces are defined as above, but replacing
$L^1$ with its natural substitute, {\it i.e.\/},
the Hardy space ~$\mathcal{H}^1$.
Similarly, in the limit case $p=\infty$ one replaces $L^\infty$ space with ${\rm BMO}$.
The classical reference for function spaces is~\cite{Tri83}.

\medskip
We will make extensive use of the Lorentz spaces $L^{p,q}$, with $1<p<\infty$ and $1\le q\le\infty$.
For completeness we recall their definition.

Let $(X, \lambda)$ be a measure space.
 Let  $f$ be  a scalar-valued $\lambda$-measurable function and
\[ \lambda_f(s) = \lambda\{ x: f(x) >s\}. \]
Then re-arrangement function $f^*$ is defined as usual by:
\[ f^*(t) = \inf\{ s: \lambda_f(s) \leq t\}. \]

By definition, for $1<p<\infty$,
\[  L^{p,q}(\R^n) = \{ f:\R^n \to \C, \;\mbox{measurable}\; \colon \|f\|_{L^{p,q}} <\infty\},\]
where
\bg
{ \|f\|_{L^{p,q}} }=\left\{
 \begin{aligned}
&\frac{q}{p} \left[\int_0^{\infty} \left( t^{\frac{1}{p} }f^*(t)\right)^q\right]^{\frac{1}{q}}, \; \mbox{if} \; q<\infty, \\ 
&\sup _{t>0}\{t^{\frac{1}{p}} f^*(t)\}, \; \mbox{if} \; q=\infty. \notag
\end{aligned} \right.
\ed
We note that it is  standard to use the above as a norm even if it does not satisfy the triangle inequality since
one can find an equivalent norm that makes the space into a Banach space.

In particular, $L^{p,\infty}$ agrees with the weak $L^{p}$  space
(or Marcinkiewicz space)
 \[ L^{p*} = \{ f:\R^n \to \C\colon f \;\mbox{measurable},\;  \|f\|_{L^{p*}} <\infty\}.\]
The quasi-norm  
\begin{equation*}
 \|f\|_{L^{p*} } = \sup_{t>0} t[ \lambda_f(t)]^{\frac{1}{p}} 
\end{equation*}
is equivalent to the norm on $L^{p,\infty}$, for $1<p<\infty$.

Our measure $\lambda$ will be chosen to be the Lebesgue measure. The Lebesgue measure of a set $A$ will be denoted by $\mes(A)$.
For  basic properties  of these spaces  useful reference are  also ~\cite{Zie}, \cite{Lem02}. It is well-known that the space $L^{p,q}$, $1<p<\infty$ and $1\leq q\leq\infty$, is the interpolated space $L^{p,q}=[L^1,L^\infty]_{1-\frac1p,q}$. Here $[\cdot
,\cdot]_{1-\frac1p,q}$ denotes the interpolated space by the real interpolation method. Using the reiteration theorem for interpolation, see \cite[Theorem 2.2]{Lem02}, one has that $L^{p,q}=[L^{p_1,q_1},L^{p_2,q_2}]_{\theta,q}$ for all $1<p_1<p_2<\infty$, $1\leq q,q_2,q_2\leq\infty$, $0<\theta<1$ and $\frac1p=\frac{1-\theta}{p_1}+\frac\theta{p_2}$. In particular, one has that $L^{p_1,q_1}\cap L^{p_2,q_2}\subset L^{p,q}$ for all  $1<p_1<p<p_2<\infty$ and $1\leq q,q_2,q_2\leq\infty$. The H\"older inequality in Lorentz spaces can be stated in the following form.
\begin{proposition}\label{holder}
Suppose that 
\begin{equation*}
1<p,p_1,p_2<\infty, \quad 1\leq q,q_1,q_2\leq\infty, \quad \frac1p=\frac1{p_1}+\frac1{p_2}\quad\text{and}\quad \frac1q=\frac1{q_1}+\frac1{q_2}.  
\end{equation*}
Then the pointwise product is a bounded bilinear operator from $L^{p_1,q_1}\times L^{p_2,q_2}$ to $L^{p,q}$, from $L^{p,q}\times L^\infty$ to $L^{p,q}$ and from $L^{p,q}\times L^{p',q'}$ to $L^1$ where $\frac1p+\frac1{p'}=1$ and $\frac1q+\frac1{q'}=1$.
\end{proposition}
The proof of this proposition can be found in \cite[Proposition 2.3]{Lem02}. The similar property for convolution is proved in \cite[Proposition 2.4]{Lem02} and reads as follows.
\begin{proposition}\label{young}
Assume  that 
\begin{equation*}
1<p,p_1,p_2<\infty, \quad 1\leq q,q_1,q_2\leq\infty, \quad 1+\frac1p=\frac1{p_1}+\frac1{p_2}\quad\text{and}\quad \frac1q=\frac1{q_1}+\frac1{q_2}.  
\end{equation*}
Then the convolution is a bounded bilinear operator from $L^{p_1,q_1}\times L^{p_2,q_2}$ to $L^{p,q}$, from $L^{p,q}\times L^{1}$ to $L^{p,q}$ and from $L^{p,q}\times L^{p',q'}$ to $L^\infty$ where $\frac1p+\frac1{p'}=1$ and $\frac1q+\frac1{q'}=1$.
\end{proposition}

\medskip

We also recall the definition of the Morrey--Campanato spaces.
In their homogeneous version,  for $1\le q\le p$, their elements are all the $L^q_{\rm loc}(\R^3)$ functions~$f$
satisfying
\begin{equation*}
\|f\|_{\mathcal{M}_{p,q}}=\sup_{x_0\in\R^3}\sup_{R>0} R^{\frac3p-\frac3q}\Bigl(\int_{|x-x_0|<R} |f(x)|^q\,dx\Bigr)^{\frac1q}<\infty
\end{equation*}
We recall that
\begin{equation}\label{mpq}
L^p=L^{p,p}=\mathcal{M}_{p,p}\subset L^{p,\infty}\subset \mathcal{M}_{p,q}, \qquad 1\le q<p<\infty,
\end{equation}
with continuous injections. 
The $\mathcal{M}_{p,q}$ spaces are of course increasing in the sense of the inclusion
as $q$ decreases. On the other hand, the $L^{p,q}$- spaces increase with~$q$.

\medskip
For $\theta\ge0$, we introduce the space $\dot E_\theta$ of all measurable functions 
(or vector field) $f$ in $\R^3$, such that
\begin{equation*}
\|f\|_{\dot E_\theta}\equiv \hbox{ess\,sup}_{x\in\R^3} \,|x|^\theta|f(x)|<\infty.
\end{equation*}

\medskip

\subsubsection{Other notations}

We denote by $\P={\rm Id}-\nabla\Delta^{-1}{\rm div}$ the Leray projector
onto the divergence-free vector field.
Notice that $\P$ is a pseudodifferential operator of order zero, which is bounded
in  $H^s_p$, $\dot H^{s}_p$ and $L^{p,q}$, for $1<p<\infty$, $1\le q\le\infty$
and $s\in\R$.
Thus, when~$f$ belongs to those spaces, the validity of an Helmholtz decomposition
$f=\P f+\nabla g$ implies that one could assume, without restriction, that $f$ is divergence-free.

However, we will  \emph{not make} this assumption in order to avoid unpleasant restrictions,
especially when working in weighted spaces
(notice that $\P$ is \emph{not} bounded in $\dot E_\theta$)
or in~$L^1$.
Indeed, it has some interest to consider integrable external forces
with non-zero mean, which prevents $\hbox{div}\,f=0$.


\section{Solutions in $L^p(\R^3)$}

The equations~\eqref{SNS} are invariant by the natural scaling $(U,p,f)\mapsto(U_\lambda,p_\lambda,f_\lambda)$
for all $\lambda>0$ and $U_\lambda=\lambda U(\lambda \cdot)$, $P_\lambda=\lambda^2P(\lambda\cdot)$
and $f_\lambda=\lambda^3f(\lambda\cdot)$.
Following a well established procedure, not only for Navier--Stokes, we consider the following program:

\medskip
(1) Existence:   first construct (rough) solutions $U$ in a scaling invariant setting,
{\it i.e.\/} in a functional space  with the same homogeneity of~$L^3$
assuming that the norm of $f$ is small in a function space (as large as possible)
with the same homogeneity of~$L^1$.

(2) Propagation: deduce from additional properties of $f$ (oscillations, localization,\dots)
additional properties for~$U$ (localization, asymptotic properties,\dots).

\medskip
We will not discuss  the propagation of the regularity since this  issue
is already  well understood
(see~\cite{GaldiBook}).
For example for, not necessarily small, external forces 
belonging to $\dot H^{-1}\cap H^s$, with $s>\frac32$,
one deduces that solutions with finite Dirichlet integral
are twice continuously differentiable and solve~\eqref{SNS} in the
classical sense.

\medskip
Concerning the first part of this program, in order to give a sense to the nonlinearity one
wants to have $U\in L^2_{\rm loc}$. As noticed in~\cite{Mey99},
the largest Banach space~$X$,  continuously included in $L^2_{\rm loc}(\R^3)$,
which is invariant under translations and such that $\|U_\lambda\|_X=\|U\|_X$,
is the  Morrey--Campanato space $\mathcal{M}_{3,2}$.
Therefore, the weakest possible smallness assumption under which one can hope
to apply the first part of the program should be
\begin{equation*}
 \|\Delta^{-1} f\|_{\mathcal{M}_{3,2}} <\ep .
\end{equation*}

However, it seems impossible to prove the existence of a solution under this type of condition.
Indeed, $U\otimes U$ would belong to~$\mathcal{M}_{\frac32,1}$, and the singular integrals
involved in equivalent formulations of~\eqref{SNS} are badly behaved in Morrey spaces of $L^{1}_{\rm loc}$
functions  (see the analysis of Taylor~\cite{Taylor92} and in particular Eq.~(3.37) of his paper).

Here the situation is less favorable than for the free \emph{non-stationary} Navier--Stokes equations,
where the existence of a global in time solution can be ensured if the initial datum of the Cauchy problem 
is small in~$\mathcal{M}_{3,2}$ (or even under more general smallness assumptions, see \cite{Lem02}).
The complication, in our case, arises from the lack of the regularizing effect of the heat kernel.

On the other hand, the above difficulty disappears in the slightly smaller spaces~$\mathcal{M}_{3,q}$.
Indeed, Kozono and Yamazaki established the following result

\begin{theorem}[See \cite{KozY95}]
\label{theoKY}
Let $2<q\le 3$.
Then there exists a positive number $\delta_q$
and a strictly monotone function~$\omega_q(\delta)$ on $[0,\delta_q]$ satisfying
$\omega_q(0)=0$, such that the following holds:
\begin{itemize}
\item
For every $f\in\mathcal{D}'(\R^3)$ there exists at most one solution~$U$ in~$\mathcal{M}_{3,q}$
satisfying $\|U\|_{\mathcal{M}_{3,q}}<\omega_q(\delta_q)$.
\item
For every tempered distribution~$f$ such that $\Delta^{-1} f\in \mathcal{M}_{3,q}$, and~$\delta=\|\Delta^{-1} f\|_{\mathcal{M}_{3,q}}<\delta_q$,
there exists a solution~$U\in \mathcal{M}_{3,q}$ of~\eqref{SNS}, such that
$\|U\|_{\mathcal{M}_{3,q}}\le \omega_q(\delta)$.
\end{itemize}
\end{theorem}

This result provides a satisfactory answer to Part~1 of the above program,
but it seems difficult to make progress in Part~2 using such functional
setting.
For example, a very strong additional condition like $f\in\mathcal{S}_0(\R^3)$
(the space of functions in the Schwartz class with vanishing moments of all order),
and $f$ small, but only  in the $\mathcal{M}_{3,q}$-norm (with~$2<q<3$), seems to imply
no interesting asymptotic properties for~$U$ (such as $U\in L^p$ with low $p$).

\medskip
On the other hand the $\mathcal{M}_{3,q}$ spaces, as $q\uparrow3$, become very close to
$L^{3,\infty}$ as can be seen from relation \eqref{mpq}.
The purpose of our first theorem is to show that one can obtain propagation results
according to Part 2 of our program, by strengthening a little the smallness assumption,
and requiring that
\begin{equation}
\label{small1}
 \|\Delta^{-1} f\|_{L^{3,\infty}}<\epunu .
\end{equation}

The continuous embedding of $L^3$ into the weak space $L^{3,\infty}$
implies that condition~\eqref{small1} will be fulfilled if, {\it i.e.\/},
$f\in \dot H^{-2}_3$ with small $\dot H^{-2}_3$-norm.
Moreover, the continuous embedding
$$\dot H^{-\frac32}\subset \dot H^{-2}_3$$
shows that  the case of forces $f\in \dot H^{-\frac32}$ with small
$\dot H^{-\frac32}$-norm is also encompassed by~\eqref{small1}.

We now state our first theorem.

\begin{theorem}
\label{theoLp}
There exists an absolute constant $\epunu >0$ with the following properties:
\begin{itemize}
\item
If $f\in \mathcal{S}'(\R^3)$ is such that $\Delta^{-1} f\in L^{3,\infty}$ and
satisfying condition~\eqref{small1},
then there exists a solution~$U\in L^{3,\infty}$ of~\eqref{SNS} such that
\begin{equation}
\label{smallU}
 \|U\|_{L^{3,\infty}}\le 2\|\Delta^{-1}\P  f\|_{L^{3,\infty}}.
\end{equation}
(The uniqueness holds in the  more general setting of Theorem~\ref{theoKY}).

\item
Let $\frac{3}{2}<p<\infty$. If~$U$ is the above solution then we have more precisely
\begin{equation*}
 U\in L^{3,\infty}\cap L^p \quad\hbox{if and only if}\quad  \P f\in \dot H^{-2}_p.
\end{equation*}
In this case
(and if $p\not=3$) $U\in L^q$ for all $q$ such that $3<q\le p$  (or $p\le q<3$).

Moreover,
$U$ belongs to $L^{3,\infty}\cap L^\infty$ (respectively, $U\in L^{3,\infty}\cap {\rm BMO}$)
 if and only if $\Delta^{-1}\P f\in L^\infty$
(respectively, $\Delta^{-1} \P f\in {\rm BMO}$).
\end{itemize}
\end{theorem}

\begin{remark}
Important examples of solutions that can be obtained through this theorem
are those corresponding to external forces $f=(f_1,f_2,f_3)$ 
with components of the form $\ep \delta$, where $\delta$ is the Dirac mass at the origin.
Notice that $f \not\in  \dot{H}^{-3/2}$. However, 
assumption~\eqref{small1} is fulfilled, because $\Delta^{-1} f(x)=\frac{\ep }{|x|}(c_1,c_2,c_3)$.

In fact, due to the invariance under rotations of~\eqref{SNS}, in this case
one can always fix a coordinate system in a way such that $f=(\ep \delta,0,0)$.
The solutions that one obtains in this way are well-known:
they are the axi-symmetric solutions (around the $x_1$ axis) discovered by Landau
sixty years ago, with ordinary differential equations methods.
These are solutions that are singular at the origin
--- in fact the components of the velocity field are homogeneous functions of degree $-1$ ---
and smooth outside zero.
They can also be seen as self-similar stationary solutions of the {\it non-stationary\/}
Navier--Stokes equations.

We refer to~\cite{CanK} for an explicit expressions and other interesting properties about these solutions and to
\cite{Sve06} (see also \cite{TianX}) for related uniqueness results.
\end{remark}

\begin{remark}
\label{remrefe}
The particular case $p=2$ is physically relevant since it corresponds
to finite energy solutions.
The conclusion  $U\in L^2$ was obtained by Bjorland and Schonbek~\cite{BjorSch07},
under 
a technical smallness assumption non invariant under scaling.
Part (2) of Theorem~\ref{theoLp} improves their result.
 Indeed the same conclusion
 can be reached under the more general conditions 
\eqref{small1} and $f\in \dot H^{-2}$.
In particular, it follows that $f\in \dot H^{-\frac32}\cap \dot H^{-2}$
with $f$ small in $\dot H^{-\frac32}$ would be enough to get $U\in L^2$.
This fact was pointed out to the first and the last author by an anonymous referee
of their paper~\cite{BjorSch07}.

Roughly speaking, for $f\in \dot H^{-\frac32}$, the additional requirement $f\in \dot H^{-2}$ (which turns out to be 
also necessary for obtaining $U\in L^2$, up
to a modification of~$f$ with an additive  potential force, which in any case
would change only the pressure of the flow),
is formally equivalent 
to the additional vanishing condition $\widehat f(\xi)=o(|\xi|^{\frac12})$ as $|\xi|\to0$. 
\end{remark}

\begin{remark}
The first conclusion of Theorem~\ref{theoLp} bears some relations
with the work of Kozono and Yamazaki~\cite{KozY98} and Yamazaki~\cite{Yam00},
where they  also obtained  existence results of (possibly non-stationary) solutions in Lorentz-spaces
and in unbounded domains.
However, the assumptions in \cite{KozY98, Yam00} on the external force reads $f=\hbox{div}\,F$,
with $F$ small in $L^{\frac32,\infty}$.
Their condition is more stringent than our condition~(\ref{small1})
because it involves more regularity  (one more derivative, or more precisely, one less anti-derivative)
on~$f$.

The first part of Theorem~\ref{theoLp} is also related
to the  work of Cannone and Karch~\cite{CanK05}.
There, the authors constructed non-stationary solutions
of Navier--Stokes in the whole space in $L^\infty_t(L^{3,\infty})$ with initial data small in $L^{3,\infty}$
and external force 
such that
$$ \sup_{t>0}\biggl\| \int_0^t e^{(t-s)\Delta}\P f(s)\,ds\biggr\|_{L^{3,\infty}}$$
is small.
With some modifications of their proofs it would be possible to deduce the first conclusion
of our theorem from their result, by considering time-independent external forces
(in this case the above condition boils down to~\eqref{small1}).
We prefer however to give a self-contained proof directly in the stationary case,
because this allows us to obtain necessary and sufficient conditions.
Moreover, none of the these papers addressed the construction of solution in $L^p$ with $p<3$.
\end{remark}

We recall a well known fixed point Lemma  for bilinear  forms that will be needed in the sequel. The proof can be found in  \cite{Can}.
 \begin{lemma} \label{Can}  
Let $X$ be a Banach space and $B : X \times X \to X$ a bilinear map. Let $\|\cdot\|_X$ denote the norm in $X$. If  for all $x_1, x_2 \in X$ one has
 \[\|B(x_1,x_2)\|_X \leq \eta \|x_1\|_X  \|x_2\|_X.\]
 Then for for all $y \in X$ satisfying $4\eta \|y\|_X <1$, the equation
 \[ x =y+ B(x,x), \]
  has a solution $x \in X$ satisfying and uniquely defined by the condition
  \[\|x\|_X \leq 2\|y\|_X. \]
\end{lemma}
\begin{remark}\label{rd}
The proof of this lemma also shows that $x=\lim\limits_{k\to\infty} x_k$ where the approximate solutions $x_k$ are defined by $x_0=y$ and $x_k=y+B(x_{k-1},x_{k-1})$.  Moreover $ \|x_k\|_X \leq 2\|y\|_X$ for all $k$.
\end{remark}

{\emph{Proof of Theorem~\ref{theoLp}}}.

We use a method of mixed bilinear estimates, inspired from~\cite{Ka84}.
Let us set
\begin{equation*}
U_0\equiv-\Delta^{-1} \P f, \qquad B(U,V)\equiv\Delta^{-1}\P\nabla\cdot (U\otimes V).
\end{equation*}
Then the system~\eqref{SNS} can be rewritten as
\begin{equation}
 \label{SNS'}
U=U_0+B(U,U)
\end{equation}
and the solutions of this equations are indeed weak solutions of~\eqref{SNS}.
This equation can be solved applying the standard fixed point method as described in Lemma \ref{Can}   in  space $L^{3,\infty}$.
We have the estimate
\begin{equation}
\label{Mey}
\|B(U,V)\|_{L^{3,\infty}}\le C_1\|U\|_{L^{3,\infty}}\|V\|_{L^{3,\infty}},
\end{equation}
for some $C_1>0$ independent on $U$ and $V$.
Note that an estimate similar to~\eqref{Mey} has been proved {\it e.g.\/} by Meyer in~\cite{Mey99}
in the case of the non-stationary Navier--Stokes equations (the bilinear operator $B$ is slightly
different in that case).

To prove~\eqref{Mey}, we only have to observe that the symbol $\widehat m(\xi)$
of the pseudo-differential operator $\Delta^{-1}\P\hbox{div}$ is a homogeneous function of degree~${-1}$,
such that $\widehat m(\xi)\in C^{\infty}(\R^3\backslash\{0\})$.
Thus, the corresponding kernel $m$ is a homogeneous function of degree $-2$,
smooth outside the origin (more precisely $m=(m_{j,h,k})_{j,h,k=1,2,3}$ and $m_{j,h,k}$ are homogeneous
functions of degree $-2$).

In particular,
$$B(U,V)=m(D)(U\otimes V) \quad \hbox{with}\quad m\in L^{\frac32,\infty}.$$
Thence,
\begin{equation}
\label{mD}
\|m(D)v\|_{L^{p_2,q_1}}\le C(p_1,q_1)\|v\|_{L^{p_1,q_1}}, \qquad 
\textstyle\frac{1}{p_2}=\frac{1}{p_1}+\frac{2}{3}-1,\quad 
\begin{cases}
1<p_1<3, \\
1\le q_1\le\infty,
\end{cases}
\end{equation}
by the Young inequality stated in Proposition \ref{young}.
Applying this to $v=U\otimes V$ and using Proposition \ref{holder} to deduce that for $U,V\in L^{3,\infty}$ one has that $v\in L^{\frac32,\infty}$,
we get estimate~\eqref{Mey} with $C_1 = C(\frac32,\infty)$.
Hence by Lemma \ref{Can} it follows that,  provided that $ 4 \|U_0\|C_1 <1$,  there exists a solution of  (\ref{SNS'}) satisfying  ~\eqref{smallU}.


\medskip
To prove Part 2,
we make use  of  approximate solutions  of  $\Phi(U)=U_0+B(U,U)$. That is we choose a sequence 
satisfying $U_k=U_0+B(U_{k-1},U_{k-1})$ and use Remark \ref{rd} to state that $U_k\to U$ in $L^{3,\infty}$ as $k\to\infty$. We show now that 
\begin{equation}
\label{mes}
 \|B(U_k,U_k)\|_{L^p}\le C(p) \|U_k\|_{L^{3,\infty}}\|U_k\|_{L^p}, \qquad \textstyle\frac{3}{2}<p<\infty,
\end{equation}
valid for some positive function $p\mapsto C(p)$, continuous on~$(\textstyle\frac{3}{2},\infty)$.
 To obtain this estimate we use the H\"older inequality given in Proposition \ref{holder} to deduce that $\|U_k \otimes U_k\|_{L^{\frac{3p}{3+p},p}}\leq C_2(p)\|U_k \|_{L^{p}}\|U_k \|_{L^{3,\infty}}$. Relation \eqref{mD} for $p_1=\frac{3p}{3+p}$ and $q_1=p$ completes the proof of \eqref{mes}. 

By Part 1, applied to the approximations $U_k$,  we know that  $\|U_k\|_{L^{3,\infty}} \leq 2\|U_0\|_{L^{3,\infty}}$. Choose  $\|U_0\|_{L^{3,\infty}}\leq c_0\epunu  $, then we get from \eqref{mes}, for $\frac{3}{2}<p<\infty$ 
\begin{equation}\label{stuff}
 \|U_{k+1}\|_{L^p}\le \|U_0\|_{L^p} +2c_0C(p)\epunu \|U_{k}\|_{L^p}.
\end{equation}

If $\P f\in\dot H^{-2}_{p}$, then $U_0\in L^p$ and so, by induction, $\|U_k\|_{L^p}<\infty$ for all~$k$.
Provided $2c_0C(p)\epunu <1$, iterating  inequality (\ref{stuff}) implies that $U_k$
is uniformly bounded in $L^p$ with respect to $k$, and hence $U\in L^p$.

However, $C(p)$ blows up as $p\to\frac{3}{2}$ or $p\to\infty$,
and we want to have a smallness assumption independent of~$p$.
To circumvent this difficulty, we replace, if necessary, the constant $\epunu $
of Part 1 of the theorem with a smaller absolute constant (still denoted~$\epunu $),
in a such way that
$2c_0\epunu <1/\sup_{2\le p\le 7}C(p)$.
Then the above argument yields the conclusion of the ``if part'' of the theorem in the case $2\le p\le 7$.
To prove the ``only if'' part one simply uses estimate~\eqref{mes} with $U_k=U$ together with {\ref{SNS'})
to get $U_0 =-\Delta(\P f) \in L^p$, hence $\P f \in\dot{H}_p^{-2}$.
\medskip

Let us now consider the case $\P f\in \dot H^{-2}_p$, $\frac{3}{2}<p<2$.
Then $U_0 \in L^p\cap L^{3,\infty}$ and by interpolation $U_0\in L^2$, so using the case $2\leq p\leq7$ we get that $U\in L^2$.
On the other hand, according to Proposition \ref{young}
the space $L^{\frac32,\infty}$ is stable under convolution
with $L^1$-functions so 
\begin{equation*}
B(U,U)=m(D)(U\otimes U)\in L^{\frac32,\infty}.
\end{equation*}
But from estimate (\ref{Mey}) we know that $B(U,U)\in L^{3,\infty}$. 
By interpolation, $ B(U,U) \in L^p$.  Combining this with equality~\eqref{SNS'} yields $U\in L^p$.  Conversely, suppose that $U\in L^p$. Since we already know that the solution $U\in L^{3,\infty}$, we deduce by interpolation that $U \in L^2$. The argument above shows that $B(U,U) \in L^p$. Hence by  (\ref{Mey}) it follows that $U_0 \in L^p$, and this, in turn is equivalent to $\P f\in\dot H^{-2}_p$.

\medskip

We now consider the case $U_0\in L^p$ with $7<p\leq\infty$ (this is equivalent to $\P f\in\dot H^{-2}_p$ if $p<\infty$).
Since $U_0\in L^{3,\infty}$, we have by interpolation that $U_0\in L^4\cap L^7$. From the previous case, we infer that  $U\in L^4\cap L^7$. By interpolation, we also have that  $U\in L^{6,2}$.
Using Proposition \ref{holder} this implies that  $U\otimes U\in L^{3,1}$, so from Proposition \ref{young} and recalling $m\in L^{\frac32,\infty}$ we get that 
\begin{equation*}
B(U,U)=m(D)(U\otimes U)\in L^\infty.
\end{equation*}
But we also know that $B(U,U)\in L^{3,\infty}$, so by interpolation $B(U,U)\in L^{p}$. From \eqref{SNS'} we conclude that $U\in L^p$.  The same argument also shows that $U\in L^p$ implies $U_0\in L^p$. 

\medskip

Finally, the BMO case follows in the same way. Indeed, the argument above shows that if $U_0$ or $U$ belong to BMO, then $B(U,U)\in L^\infty$. But $L^\infty\subset {\rm BMO}$, so $B(U,U)\in {\rm BMO}$. From relation \eqref{SNS'} we see that $U\in {\rm BMO}$ iff $U_0\in {\rm BMO}$. This completes the proof of Theorem \ref{theoLp}.
\endProof

\begin{remark}\label{remlpr}
With the same proof, one can show the following equivalent condition for the stationary solution $U$ constructed in Theorem \ref{theoLp} to belong to $L^{p,r}$. If $p\in(\frac32,\infty)$ and $r\in [1,\infty]$ then $U\in L^{p,r}$ if and only if $\Delta^{-1}\P f\in  L^{p,r}$.   
\end{remark}


\section{Pointwise behavior in $\R^3$ and asymptotic profiles}

In the previous section we dealt with forces such that $\Delta^{-1} f\in L^{3,\infty}$.
Since the typical example of a function in $L^{3,\infty}$ is $|x|^{-1}$,
it is natural to ask which supplementary properties are satisfied by the solution
when $|\Delta^{-1}f(x)|\le \ep |x|^{-1}$.
The  theorem below provides a rather complete answer.

In particular, we will obtain exact asymptotic profiles in the far field for decaying solutions of~\eqref{SNS}.
Starting with the work of Finn (see \cite{Finn65} and the references therein),
a lot is known about the spatial asymptotics of stationary solutions
in unbounded domains.
The case of the whole space that we treat in this section is of course simpler
than the case of exterior domains or aperture domains considered {\it e.g.\/} 
in~\cite{GaldiBook}. 
Nevertheless, focusing on this case allow us to put weaker (and more natural)
smallness assumptions on the force, thus providing a more
transparent presentation of the problem.

\medskip
We observe  here that, despite the unboundedness of~$\P$ in the $\dot E_\theta$
spaces,
it is fairly easy to ensure {\it e.g.\/} that $\Delta^{-1}\P f\in \dot E_1$.
Indeed, one has for example that
\begin{equation}
\label{invs}
\|\Delta^{-1}\P f\|_{\dot E_1}\le C(\|f\|_{\dot E_3}+\|f\|_{L^1}).
\end{equation}
Notice that all the norms in inequality~\eqref{invs} are invariant under scaling.
The above inequality can be proved with a simple size estimate (using that $\Delta^{-1}\P$
is a convolution operator with a kernel $\widetilde m$ satisfying $|\widetilde m(x)|\le C|x|^{-1}$).
The same conclusion $\Delta^{-1}\P f\in \dot E_1$ can be obtained also via the Fourier
transform (using classical results in~\cite{Ste93}),
assuming, {\it e.g.\/}, $f=\nabla \cdot F$ where $F$ is a two dimensional tensor
with homogeneous components of degree~$-2$, smooth outside the origin.

Let us recall the imbedding  $\dot E_1\hookrightarrow L^{3,\infty}$, thus the smallness assumption in the space $\dot E_1$ implies a smallness assumption in $L^{3,\infty}$.

The spirit of Theorem 3.1 below is close to a
previous work of the second  author (see~\cite{BraV07})
in which similar conclusions are shown for the time-dependant
Navier-Stokes equation in the whole space.

\begin{theorem}
\label{theoP}
There exists an absolute constant $\epdoi>0$ (with $\epdoi$
{\rm a priori} smaller than the constant $\epunu $
of Theorem~\ref{theoLp}) such that:
\begin{itemize}

 \item If $f\in\mathcal{S}'(\R^{3})$ is such that $\Delta^{-1} \P f\in \dot E_1$ and $\|\Delta^{-1}\P f\|_{\dot E_1}<\epdoi$,
then  the solution $U\in L^{3,\infty}$ obtained in Theorem~\ref{theoLp}
satisfies
\begin{equation*}
 \|U\|_{\dot E_1}\le 2 \|\Delta^{-1}\P f \|_{\dot E_1}.
\end{equation*}

\item
Let $0\le \theta\le 2$.
Under the additional assumption $\Delta^{-1} \P f\in \dot E_\theta$,
we have also $U\in \dot E_\theta$.

\item
In particular, if $\Delta^{-1}\P f\in \dot E_0\cap \dot E_2$, with small $\dot E_1$-norm,
then $U$ satisfies the pointwise estimate
\begin{equation*}
|U(x)|\le C(1+|x|)^{-2}.
\end{equation*}
In this case the solution $U$ has the following profile as $|x|\to\infty$:
\begin{equation}
\label{asU}
U(x)=-\Delta^{-1}\P f(x)\,+\, m(x)\colon\Bigl(\int U\otimes U\Bigr)\,+\, O\bigl(|x|^{-3}\log(|x|)\bigr),
\end{equation}
where $m=(m_{j,h,k})$ is the kernel of $\Delta^{-1}\P\hbox{div}$ 
and $m_{j,h,k}(x)$ are homogeneous functions of degree $-2$, 
$C^\infty$ outside zero. 
Furthermore,
\begin{equation}
\label{ort}
 m(x)\colon\Bigl(\int U\otimes U\Bigr)\equiv0\quad\hbox{if and only if}\quad
\exists\, c\in\R\;\hbox{s.t.} \;\int U_hU_k=c\,\delta_{h,k},
\end{equation}
for $h,k=1,2,3$, where $\delta_{h,k}=0$ or $1$ if $h\not=k$ or $h=k$.
\end{itemize}
\end{theorem}

\begin{remark}
Let us be more explicit with our notation:
by definition, for $j=1,2,3$,
$$\Bigl[m(x) \colon \int (U\otimes U)\Bigr]_j=\sum_{h,k=1}^3 m_{j,h,k}(x)\biggl(\int U_h(y)U_k(y)\;dy\biggr).$$
Moreover $m_{j,h,k}(x)=\partial_{h}M_{j,k}(x)$,
where $M_{j,k}$ is the tensor appearing in the fundamental solution of the Stokes equation.
The computation of~$M$  goes back to Lorentz (1896). See~\cite[Vol. I, p. 190]{GaldiBook}
for the explicit formula.
\end{remark}

\begin{remark}
For example, it follows from this theorem that,  if $f\in \mathcal{S}(\R^3)$
is such that  $0\not\in \hbox{supp} \widehat f$ 
and~$f$ satisfies the previous smallness assumption,
then 
\begin{equation*}
U(x)\simeq m(x)\colon\biggl(\int U\otimes U\biggr), \qquad \hbox{as $|x|\to\infty$}
\end{equation*}
provided that the right-hand side does not vanish.
Indeed, we have in this case  $\Delta^{-1}\P f\in \mathcal{S}(\R^3)$.
In particular $|U(x)|\le C(1+|x|)^{-2}$.
But the improved estimate $U(x)=o(|x|^{-2})$ as $|x|\to\infty$ holds
if and only if the flow satisfies the orthogonality relations~\eqref{ort}.
Of course, generically it is not the case.
This implies the optimality
of the restriction $\theta\le 2$ in Theorem~\ref{theoP} 
as well as the optimality of the restriction $p>\frac32$ appearing
in Theorem~\ref{theoLp}.
It is possible to relax the condition that $0\not\in\hbox{supp}\widehat f$
assuming, instead that $|\widehat f(\xi)|\le C|\xi|^k$ for a sufficiently large $k>0$.
As noticed in~\cite{BjorSch07}, this is essentially an oscillatory condition on~$f$,
describing the large time behavior of the solution of the Cauchy problem for the heat equation.
\end{remark}

\begin{remark}
Examples of (exceptional) stationary flows satisfying the orthogonality relations~\eqref{ort},
and such that $U(x)=O\bigl(|x|^{-3}\log(|x|)\bigr)$, are easily constructed by
taking $f$ satisfying the assumptions of the previous remark
and additional suitable symmetries.
An axi-symmetry condition would not be enough:
 one rather needs here polyhedral-type symmetries.
The suitable symmetries to be imposed on~$f$ can be classified exactly as done in \cite{Bra04iii}, in the
case of the non-stationary Navier--Stokes equations.
For example the two conditions 
$Rf(x)=f(Rx)$ and $Sf(x)=f(Sx)$ where $R,S$ are the orthogonal transformations
in~$\R^3$ $R\colon (x_1,x_2,x_3)\mapsto (x_2,x_3,x_1)$ and $S\colon(x_1,x_2,x_3)\mapsto(-x_1,x_2,x_3)$
are sufficient.
See~\cite{Bra04iii} for explicit examples of this type of vector fields.

On the other hand, explicit examples of solutions $U=U_f$ which {\it do not\/} satisfy the orthogonality
relations can be obtained simply by taking $f=\eta f_0$ with $\eta>0$ sufficiently small
and $f_0\in \dot H^{-2}$ satisfying the conditions of Part~2 of Theorem~\ref{theoLp} with $p=2$
(this implies that $U_{f_0}\in L^2$).
If, in addition, \emph{there is no} $c\in\R$ such that
\begin{equation*}
\int (\Delta^{-1}\P f_0)_h(\Delta^{-1}\P f_0)_k= c\,\delta_{h,k},
\end{equation*}
then $U_f$ cannot satisfy the orthogonality relations, provided $\eta>0$ is small enough.
The proof of this claim relies on an  argument that has been used in~\cite{Bra07}
in the setting of the non-stationary Navier--Stokes equations.
These observations lead us to the following theorem,
containing the announced non-existence result of generic solutions in $L^p$, $p\le \frac{3}{2}$.

\begin{theorem}
 \label{theonp}
Let $f_0=(f_1,f_2,f_3)$ be a divergence-free vector field
such that $\widehat f\in C^\infty_0(\R^3)$ and $0\not\in \hbox{supp}(\widehat{f})$.
Assume also that the matrix 
$$\biggl(\int \frac{(\widehat{f}_0)_j(\overline{\widehat{f}_0})_k}{|\xi|^4}\,d\xi\biggr)_{j,k}$$
is not a scalar multiple of the identity.
Then there exists $\eta_0>0$ such that the solution
of~\eqref{SNS} with $f=\eta f_0$ and $0<\eta\leq\eta_0$
satisfies,
\begin{equation}
\label{ulob}
c|x|^{-2}\le |U(x)|\le C|x|^{-2},
\qquad\hbox{$|x|>\!\!\!>1$},
\end{equation}
where 
$C>0$ is independent on~$x$ and $c=c(\frac{x}{|x|})$ is independent on $|x|$; 
moreover, $c(\frac{x}{|x|})>0$ on a set of positive surface measure on the unit sphere.
In particular,  $U\not\in L^p(\R^3)$ for all $1\le p\le \frac{3}{2}$.
\end{theorem}

\end{remark}

\medskip
\noindent
\emph{Proof of Theorem~\ref{theoP}. }
We already have, by Theorem~\ref{theoLp}, a solution in $L^{3,\infty}$.
To see that such solution belongs more precisely, to $\dot E_1$ we only have
to prove the estimate
\begin{equation}
 \label{biE1}
\|B(U,V)\|_{\dot E_1}\le C\|U\|_{\dot E_1}\|V\|_{\dot E_1},
\end{equation}
for some $C>0$ independent on $U$ and $V$. Indeed, an application of Lemma \ref{Can} shows the existence and the uniqueness of the solution $U$ in $\dot E_1$. This solution also belongs to $L^{3,\infty}$ since $\dot E_1\subset L^{3,\infty}$.
Of course, the re-application of the fixed point argument requires that
we replace the constant~$\epunu >0$ of Theorem~\ref{theoLp}
by a smaller one. Relation \eqref{biE1} is a particular case of the following lemma:

\begin{lemma}\label{lemmE}
Let $\theta_1,\theta_2$ be two real numbers such that $1<\theta_1+\theta_2<3$. There exists a constant $C$ such that
\begin{equation*}
\|B(U,V)\|_{\dot E_{\theta_1+\theta_2-1}}\le C\|U\|_{\dot E_{\theta_1}}\|V\|_{\dot E_{\theta_2}}.
\end{equation*}
Moreover
\begin{equation*}
\|B(U,U)\|_{\dot E_{2}}\le C(\|U\|_{\dot E_{\frac32}}^2+\|U\|_{L^2}^2).
\end{equation*}
\end{lemma}
\begin{proof}
Recall that $B(U,V)=m\ast(U\otimes V)$ with $m$ homogeneous of degree $-2$. 
Since $|m(x)|\leq C|x|^{-2}$ we can bound  
\begin{equation*}
 B(U,V)=\int
m(x-y)\colon (U\otimes V)(y)\,dy \leq C \|U\|_{\dot E_{\theta_1}}\|V\|_{\dot E_{\theta_2}} \int \frac1{|x-y|^2|y|^{\theta_1+\theta_2}}\, dy.
\end{equation*}
It is easy to show that the last integral is a function of $|x|$ homogeneous of order $1-\theta_1-\theta_2$, so it can be bounded by $C|x|^{1-\theta_1-\theta_2}$.

To show the second part, we decompose
\begin{equation*}
 B(U,V)=\Bigl(\int_{|y|\le \frac{|x|}2}+\int_{\frac{|x|}2\leq |y|}\Bigr)
m(x-y)\colon (U\otimes V)(y)\,dy
=I_1+I_2.
\end{equation*}
We have
\begin{equation*}
|I_2|\leq  C \|U\|_{\dot E_{\frac32}}^2 \int_{\frac{|x|}2\leq |y|} \frac1{|x-y|^2|y|^{3}}\, dy\leq  \frac{C}{|x|^2} \|U\|_{\dot E_{\frac32}}^2 
\end{equation*}
where we used the same scaling argument as above to deduce the last inequality. Next, we write for $I_1$
\begin{equation*}
|I_1|\leq C \int_{|y|\le \frac{|x|}2} \frac1{|x-y|^2}|U(y)|^2\, dy 
\leq \frac{C}{|x|^2}\int_{|y|\le \frac{|x|}2} |U(y)|^2\, dy 
\leq \frac{C}{|x|^2}\|U\|_{L^2}^2.
\end{equation*}
\end{proof}

Let us now prove Part~2 of Theorem~\ref{theoP}.
We have the additional information $\Delta^{-1}\mathbb{P}f\in \dot E_\theta$. We  argue as in the proof of Theorem~\ref{theoLp}.  That is  we define
$\Phi(U)=U_0+B(U,U)$ and we choose a sequence 
satisfying $U_k=\Phi(U_{k-1})$. From Lemma \ref{lemmE} we have the estimate
\begin{equation*}
 \|B(U_k,U_k)\|_{\dot E_\theta}\le C_{\theta} \|U_k\|_{\dot E_1}\|U_k\|_{\dot E_\theta}, \qquad \textstyle0<\theta<2,
\end{equation*}
for some positive function $\theta\mapsto C_\theta$, continuous on $(0,2)$.
As in Theorem~\ref{theoLp} part 2 it follows that   the sequence of approximate solutions $U_k$  remains bounded in $\dot E_\theta$, provided that $\Delta^{-1}\P f\in \dot E_\theta$, for some
$\theta\in(0,2)$, and
$$ 2 C_\theta\|\Delta^{-1}\P f\|_{\dot E_1}<1.$$
The continuity of $C_\theta$  allows to obtain the conclusion of the theorem
(with a smallness assumption independent on~$\theta$), at least for {\it e.g.\/} $\theta\in[\frac{1}{2},\frac{7}{4}]$.
We had to exclude a neighborhood of $\theta=0$ and of $\theta=2$, where $C_\theta$ blows-up.

In the case $\frac{7}{4}<\theta\le 2$, we know that $\Delta^{-1}\mathbb{P}f\in \dot E_1\cap\dot E_\theta\subset \dot E_1\cap\dot E_{\frac74} $. So, from the previous case we deduce that   the solution $U$ satisfies
$U\in \dot E_1\cap \dot E_{\frac74}\subset L^2\cap \dot E_{\frac32}$.
Using again Lemma \ref{lemmE} we infer that $B(U,U)\in \dot E_2$. But we also know that $B(U,U)\in \dot E_1$ so $B(U,U)\in \dot E_\theta$. The conclusion in the case $\frac{7}{4}<\theta\le 2$ now follows from equation~\eqref{SNS'}.

It remains to consider the case $0<\theta<\frac{1}{2}$ (the case $\theta=0$ is contained
in Theorem~\ref{theoLp}, since $\dot E_0=L^\infty$). As above, we show that $U\in \dot E_{\frac12}\cap \dot E_1$ so $U\in \dot E_{\frac{\theta+1}2}$. From Lemma \ref{lemmE} we get that $B(U,U)\in \dot E_\theta$ so $U\in \dot E_\theta$. The proof of Part 2 of Theorem~\ref{theoP} is now completed.

\medskip
Let us prove Part 3.
We will show using decay properties of $m$ and a Taylor expansion that
 for any solution such that $|U(x)|\le C(1+|x|)^{-2}$,
we have
\begin{equation}
\label{asB}
 \Delta^{-1}\P \nabla\cdot(U\otimes U)(x)=m(x)\colon\int U\otimes U +O\bigl(|x|^{-3}\log(|x|)\bigr),
\qquad \hbox{as $|x|\to\infty$}.
\end{equation}
But,
\begin{equation*}
 \begin{split}
  \Delta^{-1}\P\nabla\cdot(U\otimes U)(x)
=&\int m(x-y)\colon U\otimes U(y)\,dy\\
=& m(x)\colon\int U\otimes U \,-\,
 m(x)\colon\int_{|y|\ge |x|/2} U\otimes U \\
&+\int_{|y|\le |x|/2}[ m(x-y)-m(x)]\colon U\otimes U(y)\,dy\\
&+ \int_{|x-y|\le |x|/2} m(x-y)\colon U\otimes U(y)\,dy\\
&+\int_{|y|\ge |x|/2,\;|x-y|\ge |x|/2} m(x-y)\colon U\otimes U(y)\,dy.
 \end{split}
\end{equation*}
The only properties on the kernel $m$ that we will use are
 $|m(x)|\le C|x|^{-2}$ and $|\nabla m(x)|\le C|x|^{-3}$. 
 We need to show that all the terms on the RHS of the last inequality 
  (excepted the first one) are bounded
 by $C|x|^{-3}\log|x|$ for large $|x|$.
 This follow easily since $U \in L^2 \cap \dot{E}_2$. 
For large~$|x|$, the second, the fourth and the last term on the right-hand side
are in fact bounded by $C|x|^{-3}$.
The third term is bounded by $C|x|^{-3}\log|x|$, for large~$|x|$, as it can be checked applying the
Taylor formula to~$m$.
This implies both the asymptotic profiles~\eqref{asB} and~\eqref{asU}

\medskip
To conclude, it remains to show that the homogeneous functions
\begin{equation*}
 \sum_{h,k}m_{j,h,k}(x)\int U_hU_k, \qquad j=1,2,3,
\end{equation*}
vanish identically if and only if the matrix $\int U\otimes U$ is a scalar multiple of the identity.
We reproduce a computation similar to that in~\cite{BraV07, MiyS00}:
taking the Fourier transform, the above vanishing condition is proved
to be equivalent to
\begin{equation*}
\sum_{h,k}\widehat m_{j,h,k}(\xi)\int U_hU_k
=\sum_{h}\frac{{\rm i}\xi_h}{|\xi|^{2}}\int{U_jU_h}-\sum_{h,k}\frac{{\rm i}\xi_j\xi_h\xi_k}{|\xi|^4}\int U_hU_k=0,
\qquad\hbox{for a.e. $\xi\in\R^3$}.
\end{equation*}
The conclusion is now obvious.

\endProof

We end this section with the proof of Theorem \ref{theonp}.

\begin{proof}[Proof of Theorem \ref{theonp}]
We start by choosing $\eta_0$ sufficiently small such that 
$$\eta_0\|\Delta^{-1}\mathbb{P}f_0\|_{L^{3,\infty}}\leq \epunu ,$$
where $\epunu $ is the smallness constant from Theorem \ref{theoLp}. According to Theorem \ref{theoLp},
for $0<\eta\leq\eta_0$ there exists a unique solution $U\in L^{3,\infty}\cap L^2$ of \eqref{SNS} with $f=\eta f_0$
 such that $\|U\|_{L^{3,\infty}}\leq  2\eta \|\Delta^{-1}\mathbb{P}f_0\|_{L^{3,\infty}}$.
It suffices to show that the orthogonality relations \eqref{ort} does not hold true for $U$.  

Let $W_0= -\Delta^{-1}\mathbb{P}f_0$ and $U_0=\eta W_0$. The hypothesis implies that the matrix $\int W_0\otimes W_0$ is not a scalar multiple of the identity. This means that there exists $j\neq k$ such that either $\int W_{0}^j W_{0}^k\neq0$ or $\int |W_{0}^j|^2\neq \int |W_{0}^k|^2$, where $W_0^j$ denotes the $j$-th component of $W_0$. We will suppose that $\int W_{0}^j W_{0}^k\neq0$, the other case being entirely similar.

We have
\begin{multline}\label{uuz}
\Bigl|\int U^j U^k-\int U_{0}^j U_{0}^k\Bigr|
=\Bigl|\int (U^j-U_0^j) U^k+\int U_{0}^j (U^k-U_{0}^k)\Bigr| \\
\leq \|U-U_0\|_{L^2}(\|U\|_{L^2}+ \|U_0\|_{L^2})
\end{multline}

From \eqref{SNS'} and \eqref{mes} with $p=2$ and $U_k$ replaced by $U$ we deduce that
\begin{equation}\label{uuzz}
\|U-U_0\|_{L^2}=\|B(U,U)\|_{L^2}\leq C(2) \|U\|_{L^2}\|U\|_{L^{3,\infty}}
\leq 2C(2)\eta \|W_0\|_{L^{3,\infty}} \|U\|_{L^2}.
\end{equation}
Therefore
\begin{equation*}
\|U\|_{L^2} \leq  \|U_0\|_{L^2}+\|U-U_0\|_{L^2}
\leq  \eta\|W_0\|_{L^2}+2C(2)\eta_0 \|W_0\|_{L^{3,\infty}} \|U\|_{L^2}.
\end{equation*}
If we further strengthen the smallness assumption on $\eta_0$ by
\begin{equation*}
  \eta_0\leq \frac1{4C(2) \|W_0\|_{L^{3,\infty}}}
\end{equation*}
we get that
\begin{equation*}
\|U\|_{L^2} \leq  2\eta\|W_0\|_{L^2}. 
\end{equation*}
Relation \eqref{uuzz} combined with the previous estimate implies that 
\begin{equation*}
\|U-U_0\|_{L^2}\leq 4C(2)\eta^2 \|W_0\|_{L^{3,\infty}} \|W_0\|_{L^2}.  
\end{equation*}
Using the two previous bounds in \eqref{uuz} implies that
\begin{equation*}
\Bigl|\int U^j U^k-\eta^2\int W_{0}^j W_{0}^k\Bigr|
\leq 12C(2)\eta^3 \|W_0\|_{L^{3,\infty}} \|W_0\|^2_{L^2}.  
\end{equation*}

Finally
\begin{multline*}
\Bigl|\int U^j U^k \Bigr|\geq \eta^2 \Bigl|\int W_{0}^j W_{0}^k\Bigr| - \Bigl|\int U^j U^k-\eta^2\int W_{0}^j W_{0}^k\Bigr|\\
\geq \eta^2 \Bigl|\int W_{0}^j W_{0}^k\Bigr| - 12C(2)\eta^3 \|W_0\|_{L^{3,\infty}} \|W_0\|^2_{L^2}>0
\end{multline*}
if we further assume that 
\begin{equation*}
  \eta_0\leq \frac{\bigl|\int W_{0}^j W_{0}^k\bigr|}{24C(2) \|W_0\|_{L^{3,\infty}}\|W_0\|^2_{L^2}}.
\end{equation*}
\end{proof}

\section{Stability of the Stationary Solutions}

Consider now a mild formulation of the Navier-Stokes equations with time independent forcing function $f$
satisfying, as usual, to a smallness condition as in~\eqref{small1},
\begin{equation}
u(t)= e^{t\Delta}u_0 +\int_0^te^{(t-s)\Delta}\mathbb{P}f \,ds -\int_0^te^{(t-s)\Delta}\mathbb{P}\nabla\cdot (u\otimes u)(s) \,ds.
\label{mildNS:PDE}
\end{equation}
The two main goals of this section are the following. 
First we want to establish conditions on $u_0$ to ensure that
the above system has a solution $u\in L^\infty(\mathbb{R}_+,L^{p,\infty})$. 
Next, we want to find the largest possible class of solutions~$u$ to~\eqref{mildNS:PDE} for which we can say that~$u(t)$
converges to   the steady solution~$U$ given by Theorem \ref{theoLp} corresponding to the same force~$f$.
 This class will be general enough to include non-stationary solutions in $L^{3,\infty}$ with large initial data.
We will show in particular that  \textit{a priori} global solutions, verifying a mild regularity condition
but initially large in  $L^{3,\infty}$, become small in  $L^{3,\infty}$ after some time. Only the singularity at infinity of the initial velocity needs to be small in some sense which is made rigorous in  \eqref{smallcondl3}. For example, we allow an initial velocity $u_0$ bounded by $C/|x|$ everywhere and bounded by $\varepsilon/|x|$ for large $x$, with $C$ arbitrary and $\varepsilon$ small. 

 We recall that \textit{a priori} large non-stationary solutions
 in  $\dot B^{-1+\frac3p}_{p,q}$ and  $VMO^{-1}$ of the Navier--Stokes equations without forcing
are known to converge to zero in these spaces (see \cite{GIP,ADT}). 
However, in our case, convergence to zero will not necessarily
hold true for $\|u(t)-U\|_{L^{3,\infty}}$, due to the fact that the smooth function in $\mathcal{S}(\R^3)$ are not dense in $L^{3,\infty}$.
Thus, only weaker convergence results should be expected.

Theorem \ref{forcedNS:theorem} collects our results on the stability of small solutions~$u$, extending, 
for flows in $\R^3$ with time independent forcing term,
those of~\cite{Bar99, CanK05, KozY98, Yam00} to the case $\frac{3}{2}<p<3$, and providing some additional information also for $p>3$. Theorem \ref{largesol} contains the convergence result of large solutions~$u$ to small stationary solutions~$U$.
Its proof relies on some energy estimates inspired by \cite{GIP,ADT} and on the results
on the stability of small solutions prepared in Theorem \ref{forcedNS:theorem}.

To begin we first recall a lemma which will be useful for estimating the integral terms on the RHS of (\ref{mildNS:PDE}) in
 $L^{p,\infty}$ spaces. We notice that for the case $p=3,\, q=3/2$ the Lemma below was obtained in several papers, among them the first seems to be in Yamazaki's paper  \cite{Yam00}. 
Variants of this lemma can also be found in~\cite{Mey99}, in a slightly less general form, and in~\cite{Lem02}.

\begin{lemma}\label{heatforce:lemma}
Given any $p\in(\frac{3}{2},\infty)$ let $q=\frac{3p}{p+3}$.  For $0\le \sigma<t$, the operator
\begin{equation}
 \tilde{L}_\sigma (\phi)(t) = \int_{\sigma}^t e^{(t-s)\Delta}\mathbb{P} \nabla \cdot \phi(s)\, ds\notag
\end{equation}
satisfies
\begin{equation}
\|\tilde{L}_\sigma (\phi)(t)\|_{L^{p,\infty}}\leq C(p)\sup_{0<s<t}\|\phi(s)\|_{L^{q,\infty}}\label{heatforce:bound}
\end{equation}
where $C(p)$ denotes a constant independent of $\sigma$.
\end{lemma}
\begin{proof}
Let $F(t)$ be the kernel of the operator $e^{t\Delta}\mathbb{P}\hbox{div}$.
First recall the rescaling relation 
$$F(x,t)=t^{-2}F(x/\sqrt t,1)$$
and that $F(\cdot,1)\in L^1\cap L^\infty$.

We consider separately the following two pieces.
\begin{equation*}
A_1 = \int_{t-\lambda^*}^t F(t-s)\ast \phi(s)\, ds
\quad\text{and}\quad
A_2 = \int^{t-\lambda^*}_\sigma F(t-s)\ast \phi(s)\, ds.
\end{equation*}
The idea of the estimate is to find, given any fixed $\lambda$, a $\lambda^*$ so that $|A_2|<\lambda/2$.
 With this choice of $\lambda^*$ we can estimate the Lebesgue measure of the set $\{x:\tilde{L}_\sigma (\phi)|>\lambda\}$
 in terms of $A_1$ only.
 In that direction we establish two preliminary estimates.  The first is a an application of Young's inequality stated in Proposition \ref{young}:
\begin{equation}
\|A_2\|_{L^\infty} \leq C\int_\sigma^{t-\lambda^*}\|F(t-s)\|_{L^{\alpha,1}}\|\phi\|_{L^{q,\infty}}\,ds\notag
\end{equation}
Here, $\alpha=\frac{3p}{2p-3}$.  The estimate $\|F(t-s)\|_{L^{\alpha,1}}\leq C(t-s)^{-1-\frac{3}{2p}}$ 
(that follows from the rescaling properties of~$F$)
implies
\begin{equation}
\|A_2\|_{L^\infty}\leq C(p)(\lambda^*)^{-\frac{3}{2p}}\|\phi\|_{X^{\sigma,t}_q}\label{A2bound:hflemma}.
\end{equation}
Here we have introduced the notation $X^{\sigma,t}_q=L^{\infty}((\sigma,t),L^{q,\infty})$.  Similarly, $\|F(t-s)\|_{L^1}\leq (t-s)^{-\frac{1}{2}}$ and
\begin{equation}
 \|A_1\|_{L^{q,\infty}}\leq \int_{t-\lambda^*}^t\|F(t-s)\|_{L^1}\|\phi\|_{L^{q,\infty}}\,ds\leq (\lambda^*)^{\frac{1}{2}}\|\phi\|_{X^{\sigma,t}_q}.
\label{A1bound:hflemma}
\end{equation}

We proceed with the bound for $\|\tilde{L}_\sigma (\phi)\|_{L^{p,\infty}}$.  Using the definition of the norm and the triangle inequality,
\begin{equation}
\|\tilde{L}_\sigma (\phi)(t)\|_{L^{p,\infty}}\leq \sup_{\lambda>0} \lambda \mes\{x:|A_1| +|A_2|>\lambda\} ^{\frac{1}{p}}\notag
\end{equation}
For each $\lambda>0$ we may choose $\lambda^*$ such that the RHS of (\ref{A2bound:hflemma}) is equal to $\lambda/2$.  With this choice of $\lambda^*$,
\begin{equation}
\lambda\mes\{x:|A_1|+|A_2|>\lambda \}^{\frac1p} \leq \lambda\mes\{x:|A_1|>\lambda/2\}^{\frac1p}\notag
\end{equation}
Also, using (\ref{A1bound:hflemma}):
\begin{equation*}
\lambda\mes\{x:|A_1|>\lambda/2\}^{\frac{1}{p}} \leq \lambda^{1-\frac{q}{p}}\|A_1\|_{L^{q,\infty}}^{\frac{q}{p}}\leq  C\|\phi\|_{X^{\sigma,t}_q.}
\end{equation*}
Taking the supremum over all $\lambda>0$ establishes (\ref{heatforce:bound}).

\end{proof}

\medskip
The following lemma concerns the large time behavior in $L^{3,\infty}$ of solutions of the heat equation.
It will provide a better understanding of the statements of our two next theorems.

\begin{lemma}\label{decomp}
Let $f\in L^{3,\infty}$.
\begin{itemize}
\item Let  $\ep >0$ be arbitrary. Then $f$ can be decomposed as $f=f_1+f_2$ with $f_1\in L^2$ and $\|f_2\|_{L^{3,\infty}}<\ep$ if and only if $\limsup\limits_{R\to0}R\mes\{|f|>R\}^{\frac13}<\ep.$
\item If $\lim\limits_{R\to0}R\mes\{|f|>R\}^{\frac13}=0$ then $e^{t\Delta}f\to0$  in  $L^{3,\infty}$ as $t\to\infty$.
\item There exists some $g\in L^{3,\infty}$ such that $e^{t\Delta}g\to0$  in  $L^{3,\infty}$ as $t\to\infty$ and such that $\limsup\limits_{R\to0}R\mes\{|g|>R\}^{\frac13}\neq0$.
\end{itemize}
\end{lemma}
\begin{proof}
Assume first that  $f=f_1+f_2$ with $f_1\in L^2$ and $\|f_2\|_{L^{3,\infty}}\leq\ep$. We estimate
\begin{equation*}
\mes\{|f_1|>R\}\leq \frac1{R^2}\int_{\R^3}|f_1|^2  
\end{equation*}
so that $\limsup\limits_{R\to0}R\mes\{|f_1|>R\}^{\frac13}=0.$ We also have that 
\begin{equation*}
\limsup\limits_{R\to0}R\mes\{|f_2|>R\}^{\frac13}\leq  \sup_{R>0}R\mes\{|f_2|>R\}^{\frac13} =\|f_2\|_{L^{3,\infty}}<\ep.
\end{equation*}
Let $\delta\in(0,1)$. Since $\{|f|>R\}\subset \{|f_1|>\delta R\}\cup \{|f_2|>(1-\delta) R\}$ we have that
\begin{align*}
 \limsup_{R\to0}R\mes\{|f|>R\}^{\frac13} 
&\leq \limsup_{R\to0}R\bigl(\mes\{|f_1|>\delta R\}+\mes\{|f_2|>(1-\delta)R\}\bigr)^{\frac13} \\
&\leq \limsup_{R\to0}R\mes\{|f_1|>\delta R\}^{\frac13} + \limsup_{R\to0}R\mes\{|f_2|>(1-\delta)R\}^{\frac13} \\
&=\frac1\delta \limsup_{R\to0}R\mes\{|f_1|>R\}^{\frac13} + \frac1{1-\delta}\limsup_{R\to0}R\mes\{|f_2|>R\}^{\frac13} \\
&\leq \frac{1}{1-\delta}\|f_2\|_{L^{3,\infty}}
\end{align*}
Leting $\delta\to0$ implies that $\limsup\limits_{R\to0}R\mes\{|f|>R\}^{\frac13}\leq\|f_2\|_{L^{3,\infty}}<\ep.$

Conversely, assume that $\limsup\limits_{R\to0}R\mes\{|f|>R\}^{\frac13}<\ep.$ There exists $R_\ep$ such that 
\begin{equation*}
\sup _{0<R<R_\ep}R\mes\{|f|>R\}^{\frac13} <\ep.
\end{equation*}
We set $f_1=f\chi_{\{|f|>R_\ep\}}$ and $f_2=f\chi_{\{|f|\leq R_\ep\}}$ where $\chi$ denotes the characteristic function. Clearly $|f_2|\leq R_\ep$ and $|f_2|\leq |f|$ so that
\begin{equation*}
\|f_2\|_{L^{3,\infty}}= \sup _{0<R<R_\ep}R\mes\{|f_2|>R\}^{\frac13} 
\leq \sup _{0<R<R_\ep}R\mes\{|f|>R\}^{\frac13} <\ep. 
\end{equation*}
It remains to show that $f_1\in L^2(\R^3)$. Let $N_\ep\in\mathbb{Z}$ be such that $R_\ep>2^{N_\ep}$. Then
\begin{equation*}
\{|f|>R_\ep\}\subset \bigcup_{n\geq N_\ep} \{2^n<|f|\leq 2^{n+1}\}
\end{equation*}
so
\begin{multline*}
\int_{\R^3}|f_1|^2=\int_{\{|f|>R_\ep\}}|f|^2 
\leq \sum_{n=N_\ep}^\infty\int_{\{2^n<|f|\leq 2^{n+1}\}}  |f|^2 
\leq \sum_{n=N_\ep}^\infty 4^{n+1}\mes{\{2^n<|f|\}} \\  
\leq \sum_{n=N_\ep}^\infty \frac4{2^n} \|f\|_{L^{3,\infty}}^3<\infty.   
\end{multline*}
This shows the first part of the lemma.

\medskip

Assume now that $\lim\limits_{R\to0}R\mes\{|f|>R\}^{\frac13}=0$ and let $\ep>0$ be arbitrary. Using the first part we decompose  $f=f_1+f_2$ with $f_1\in L^2$ and $\|f_2\|_{L^{3,\infty}}<\ep$. The standard decay estimates for the heat equation implies that $\|e^{t\Delta}f_1\|_{L^{3,\infty}}<Ct^{-\frac14}\|f_1\|_{L^2}\to 0 $ as $t\to\infty$. Moreover, $\|e^{t\Delta}f_2\|_{L^{3,\infty}}\leq \|f_2\|_{L^{3,\infty}} <\ep$. We infer that
\begin{equation*}
\limsup _{t\to\infty} \|e^{t\Delta}f\|_{L^{3,\infty}}\leq \limsup _{t\to\infty} (\|e^{t\Delta}f_1\|_{L^{3,\infty}}+\|e^{t\Delta}f_2\|_{L^{3,\infty}})\leq\ep.  
\end{equation*}
Letting $\ep\to0$ yields $\limsup\limits _{t\to\infty} \|e^{t\Delta}f\|_{L^{3,\infty}}=0$, as required.

\medskip

To prove the third part of the lemma, we choose
\begin{equation*}
  g(x)=\frac{e^{i|x|^2}}{\langle x\rangle}, \qquad \langle x\rangle=(1+|x|^2)^{\frac12}.
\end{equation*}
It is a straightforward calculation to check that 
\begin{equation*}
\limsup\limits_{R\to0}R\mes\{|g|>R\}^{\frac13}=\bigl(\frac{4\pi}3\bigr)^{\frac13}.  
\end{equation*}
On the other hand, we will show that $e^{\frac14\Delta} g\in L^2$ which by the decay estimates for the heat equation implies that $\|e^{t\Delta}g\|_{L^{3,\infty}}<C(t-\frac14)^{-\frac14}\|e^{\frac14\Delta} g\|_{L^2}\to 0 $ as $t\to\infty$. 

Since the kernel of the operator $e^{\frac14\Delta}$ is $\pi^{-\frac32}e^{-|x|^2}$ one has that
\begin{align*}
e^{\frac14\Delta}g(x)=\pi^{-\frac32}\int_{\R^3}\frac{e^{i|x-y|^2}}{\langle x-y \rangle}e^{-|y|^2}\,dy
&= \pi^{-\frac32}e^{i|x|^2}\int_{\R^3} e^{-2ix\cdot y}\, \frac{e^{(i-1)|y|^2}}{\langle x-y \rangle}\,dy \\
&=\pi^{-\frac32}\frac{e^{i|x|^2}}{\langle x \rangle^2}\int_{\R^3} (1-\frac14\Delta_y)e^{-2ix\cdot y}\, \frac{e^{(i-1)|y|^2}}{\langle x-y \rangle}\,dy\\
&=\pi^{-\frac32}\frac{e^{i|x|^2}}{\langle x \rangle^2}\int_{\R^3} e^{-2ix\cdot y}\, (1-\frac14\Delta_y)\Bigr[\frac{e^{(i-1)|y|^2}}{\langle x-y \rangle}\Bigl]\,dy
\end{align*}
The integral in the last term is bounded uniformly with respect to $x$. Indeed, all derivatives of $e^{(i-1)|y|^2} $ are integrable and all derivatives of $\frac{1}{\langle x-y \rangle} $ are uniformly bounded in $x$ and $y$. We deduce that $|e^{\frac14\Delta}g(x)|\leq C\langle x \rangle^{-2} $ which implies that $e^{\frac14\Delta} g\in L^2$. This completes the proof of the lemma.
\end{proof}

We state now our stability result for small solutions.

\begin{theorem}\label{forcedNS:theorem}
There exists an absolute constant $\eptrei >0$ 
with the following properties:
\begin{itemize}
\item
If $f,u_0\in \mathcal{S}'(\mathbb{R}^3)$ are such that
\begin{equation}
\|\Delta^{-1}\mathbb{P}f\|_{L^{3,\infty}}+\|u_0\|_{L^{3,\infty}}<\eptrei  \label{small:assumption}
\end{equation}
then there is a unique solution $u\in L^\infty(\mathbb{R}_+,L^{3,\infty})$ of (\ref{mildNS:PDE}),
weakly continuous with respect to $t\in [0,\infty)$, satisfying
\begin{equation}
\sup_{s>0}\|u(s)\|_{L^{3,\infty}}\leq 2\|u_0\|_{L^{3,\infty}}+4\|\Delta^{-1}\mathbb{P}f\|_{L^{3,\infty}}.
\label{twicedata:bound}
\end{equation}
\item
Let $p\in (\frac{3}{2},\infty)$ and suppose in addition to \eqref{small:assumption} that $u_0\in L^{p,\infty}$.  If $u$ is the above solution then,
\begin{equation*}
u\in L^\infty(\mathbb{R}_+,L^{p,\infty}) \ \ \hbox{if  and  only  if} \ \ \Delta^{-1}\mathbb{P}f \in L^{p,\infty}.
\end{equation*}
\item
Let  $p\in (\frac{3}{2},\infty)$, and $q>\min\{3,p\}$. Suppose in addition to \eqref{small:assumption}
that $u_0\in L^{p,\infty}$ and $\Delta^{-1}\mathbb{P}f \in L^{p,\infty}$.
Let also $U\in L^{3,\infty}\cap L^{p,\infty}$ be the unique stationary solution given by  Theorem \ref{theoLp} (see also Remark \ref{remlpr}).
(We assume here that $\varepsilon_3\le\varepsilon_1$, the constant introduced in Theorem~\ref{theoLp}).
\begin{enumerate}
\item[(i)]
There is a function $\varepsilon(q)>0$ such that if $\varepsilon_3<\varepsilon(q)$ then,
for some constant $C>0$,
\begin{equation}
 \label{strong-q3p}
\|u(t)-U\|_{L^q} \le Ct^{-\frac{3}{2}(\frac{1}{\min(3,p)}-\frac{1}{q})},
\qquad \forall \min\{p,3\}<q<\infty.
\end{equation}
In particular, $u(t)-U\to 0$ in $L^{q}$ as $t\to\infty$ for all $q>\min\{3,p\}$.
\item[(ii)]
If $\frac{3}{2}<p\le 3$, 
then $u(t)\rightharpoonup U$ weakly in $L^{p,\infty}$ as $t\to\infty$.
Moreover, $u(t)\to U$ strongly in $L^{p,\infty}$ if and only if $e^{t\Delta}(u_0-U)\to 0$ in $L^{p,\infty}$.
\item[(iii)]
If $\frac{3}{2}<p<3$, then the conclusion of the previous item can be strengthened as follows:
\begin{equation}
\label{lhsa}
\|u(t)-U-e^{t\Delta}(u_0-U)\|_{L^q}\le Ct^{\frac12+\frac3{2q}-\frac3p}
\end{equation}
for all $\frac{3p}{6-p}\leq q\leq p$ and for some constant $C>0$ independent of~$t$.
\end{enumerate}
In particular, $u(t)-U\to0$ in $L^q$ if and only if $e^{t\Delta}(u_0-U)\to0$ in $L^q$ as $t\to\infty$,
for all $\frac{3p}{6-p}<q\leq p$.
\end{itemize}
\end{theorem}

\medskip
Notice that in~\eqref{strong-q3p} neither $u(t)$ nor $U$ belong in general to $L^q$.
Similarly, the terms appearing in the LHS of~\eqref{lhsa} in general do not belong, separately, to $L^q$.
In other words, the difference $u(t)-U$ is better behaved than the solutions themselves. 

\begin{remark}
In the particular case $p=q=2$, the preceding theorem contains an interesting variant
 of the stability result for finite-energy solutions obtained in \cite{BjorSch07} with a different method.
Indeed, consider a stationary solution $U\in L^2\cap L^{3,\infty}$ and a perturbation $w_0\in L^2\cap L^{3,\infty}$.
According to conclusion (iii), the solution $u$ of the non-stationary Navier--Stokes equations starting from~$u_0=U-w_0$
satisfies, under the above smallness assumptions, $u(t)\to U$ in $L^2$ as $t\to\infty$
(we use here that $e^{t\Delta}w_0\to 0$ in $L^2$).
Explicit convergence rates can be given, e.g., if the perturbation belongs to additional function spaces.
For instance, when $w_0\in L^{\frac{3}{2},\infty}\cap L^{3\infty}$, then
\begin{equation*}
 \|u(t)-U\|_2\le Ct^{-1/4}, \qquad \hbox{as $t\to\infty$}.
\end{equation*}
\end{remark}

\medskip
\begin{remark}
Let us  present some further immediate consequences of this theorem.
If the perturbation satisfies $w_0\in L^3$, then $e^{t\Delta}w_0\to0$ in $L^3$ and so in $L^{3,\infty}$ as $t\to\infty$.
This in turn implies, by (ii),
\begin{equation*} 
u(t)\to U \quad\hbox{in $L^{3,\infty}$ as $t\to \infty$}.
\end{equation*}
More generally, according to the second part of Lemma~\ref{decomp}, such conclusion remains valid when
$\lim\limits_{R\to0}R\mes\{|w_0|>R\}^{\frac13}=0$.
However, notice that neither $w_0\in L^{3,\infty}$ is sufficient nor
$\lim\limits_{R\to0}R\mes\{|w_0|>R\}^{\frac13}=0$ is necessary to ensure this result.

In the same way, in the case $\frac{3}{2}<p<3$, the condition 
$\lim\limits_{R\to0}R\mes\{|w_0|>R\}^{\frac1p}=0$ implies that
$u(t)\to U\in L^{p,\infty}$ as $t\to\infty$. But in this case the stronger condition $w_0\in L^p$ would imply also,
by (iii),  the stronger conclusion $u(t)\to U$ in $L^p$.
\end{remark}
\begin{remark} The proof  of Theorem \ref{forcedNS:theorem} below will show that Equation \eqref{lhsa} holds true for the wider range $\max(1,\frac p2)\leq q < \frac{3p}{3-p}$. We did not state the full range for $q$ because the most interesting case is $q\leq p$ and also because it would require showing that in the case $p<3$, the statement (i) is true with a constant $\varepsilon(q)$ independent of $q$. 
This additional fact is easy to prove with a recursive argument, but since it is not really necessary we prefer to skip it.   
\end{remark}

{\emph{Proof of Theorem~\ref{forcedNS:theorem}}}.
We estimate the forcing term in (\ref{mildNS:PDE}) by integrating the heat kernel in time then relying
on a fixed point argument making use of Lemma \ref{heatforce:lemma}.  The relation 
\begin{equation}
\int_0^te^{(t-s)\Delta}\,ds=e^{t\Delta}\Delta^{-1}-\Delta^{-1}\label{forcingterm}
\end{equation}
 that follows since both operators have the same symbols, gives
\begin{equation}\label{canka}
\bigl\|\int_0^te^{(t-s)\Delta}\mathbb{P}f \,ds\bigr\|_{L^{3,\infty}} \leq 2\|\Delta^{-1}\mathbb{P}f\|_{L^{3,\infty}}.
\end{equation}
We used above that $e^{t\Delta}$ is a convolution operator with a function of norm $L^1$ equal to 1. Given \eqref{canka}, the first part of this theorem follows from the work of Cannone and Karch~\cite{CanK05}. But the proof takes only a few lines, so we give it for the sake of the completeness.

Using again that the kernel of  $e^{t\Delta}$ is of $L^1$ norm equal to 1 we deduce that $\|e^{t\Delta}u_0\|_{L^{3,\infty}}\leq \|u_0\|_{L^{3,\infty}} $. Therefore, if we denote  $\tilde{u}_0 =e^{t\Delta}u_0 +\int_0^te^{(t-s)\Delta}\mathbb{P}f \,ds$ one has that
\begin{equation*}
\|\tilde{u}_0\|_{L^\infty(\R_+;L^{3,\infty})}\leq \|u_0\|_{L^{3,\infty}} +  2\|\Delta^{-1}\mathbb{P}f\|_{L^{3,\infty}}. 
\end{equation*}

To apply the fixed point argument we  introduce the notation $\tilde{B}(u,v)=\tilde{L}_0 (u\otimes v)$ and rewrite (\ref{mildNS:PDE}) as
\begin{equation}
u=\tilde{u}_0-\tilde{B}(u,u)\label{mildNS:PDE2}
\end{equation}
The bound (\ref{heatforce:bound}), with $p=3$, and hence $q =\frac{3}{2}$ combined with the H\"{o}lder inequality from Proposition \ref{holder} yields
\begin{equation}
\|\tilde{B}(u,v)(t)\|_{L^{3,\infty}}
\leq C\Bigl(\sup_{s>0}C\|u(s)\|_{L^{3,\infty}}\Bigr)\Bigl(\sup_{s>0}\|v(s)\|_{L^{3,\infty}}\Bigr)\notag
\end{equation}

We apply  this estimate combined with the fixed point argument given in Lemma  \ref{Can}
 to the operator $\tilde{\Phi}(u)=\tilde{u}_0-\tilde{B}(u,u)$ in the space $L^\infty(\mathbb{R}_+,L^{3,\infty})$.
 This approach yields the existence of a unique solution $u\in L^\infty(\mathbb{R}_+,L^{3,\infty})$ 
provided  $4C\|\tilde{u}_0\|_{L^\infty(\R_+; L^{3,\infty})} <1$. Lemma
  \ref{Can} also insures that the solution satisfies $\|u\|_{L^\infty(\R_+; L^{3,\infty})} \leq 2 \|\tilde{u}_0\|_{L^\infty(\R_+; L^{3,\infty})}$,
 establishing part 1 of the Theorem.

\bigskip

To prove the second part of the theorem
we  establish first the cases $p\in [2,7]$ with
a fixed point argument then treat the other cases with an interpolation argument.  First, combine (\ref{heatforce:bound}) with the H\"{o}lder inequality to establish
\begin{equation}
\|\tilde{B}(u,u)(t)\|_{L^{p,\infty}}\leq C(p)
\Bigl(\sup_{s>0}\|u(s)\|_{L^{p,\infty}}\Bigr)\Bigl(\sup_{s>0}\|u(s)\|_{L^{3,\infty}}\Bigr)
\label{bilinearp:bound}
\end{equation}
Let $\tilde{C}$ be the maximum value of the constant in the above equation for $p\in [2,7]$, we require $8\eptrei  <1/\tilde{C}$.  Considering again the sequence of approximate solutions $(u_i)$
constructed in the usual way, and making use of (\ref{twicedata:bound}) we see
\begin{equation}
\|u_{i+1}(t)\|_{L^{p,\infty}}\leq \sup_{s>0}\|\tilde{u}_0\|_{L^{p,\infty}}+\frac{1}{2}\sup_{s>0}\|u_i(s)\|_{L^{p,\infty}}\notag
\end{equation}
From this estimate the ``if'' statement in the second claim follows for $p\in [2,7]$.

If $p\in (\frac{3}{2}, 2)$, through interpolation we find that for all $r\in(2,3)$ we have that $\tilde{u}_0\in L^\infty(\R_+;L^{r,\infty})$ and therefore $u\in L^\infty(\R_+;L^{r,\infty})$.
Appealing to (\ref{heatforce:bound}) and again combining it with the H\"older inequality we see
\begin{equation}
\|\tilde{B}(u,u)(t)\|_{L^{p,\infty}}\leq C\sup_{t>0}\|u(s)\|_{L^{r,\infty}}^2\label{step1}
\end{equation}
where $r=\frac{6p}{p+3}\in (2,3)$, hence the right hand side is bounded.  Combining this estimate with (\ref{mildNS:PDE2}) is enough to prove the ``if'' statement in the case $p\in(\frac{3}{2},2)$.  If $p\in (7,\infty)$ we again interpolate to get
$u\in L^\infty(\R_+;L^{r,\infty})$ for all $r\in(3,6)$. Choosing again $r=\frac{6p}{p+3}\in (3,6)$ in (\ref{step1}) finishes the ``if'' statement in the second claim.
To establish the ``only if'' part of the claim combine (\ref{bilinearp:bound}) with (\ref{mildNS:PDE2}) and notice that the RHS of (\ref{forcingterm}) tends to $-\Delta^{-1}$ as $t\rightarrow\infty$.
The weak continuity $u(t)\to u(t')$ for $t\to t'$ and $t'\in[0,\infty)$ 
 (the continuity is actually in the strong topology of $L^{3,\infty}$ for $t'\in (0,\infty)$\,)
is proved as in~\cite{Mey99}.

\bigskip

It remains to prove the third part of the theorem, the stability results for stationary solutions.
We begin with Claim (i).
Let $q>\min\{3,p\}$.
It is worth noticing that for $q>3$, a stability result  in the $L^{q,\infty}$-norm, as well a decay estimate of the form
$\|u(t)-v(t)\|_{L^q} \le Ct^{-\frac{3}{2}(\frac{1}{3}-\frac{1}{q})}$
was stated  in~\cite[Proposition 4.3]{CanK05}. However, it seems  that the argument briefly sketched in~\cite{CanK05}
cannot be directly applied to the case where the second solution $v(t)$ is stationary, because a
non-obvious generalization of Lemma~\ref{heatforce:lemma} would be needed.
Therefore, we provide a detailed proof of estimate~\eqref{strong-q3p}.

Let $w=U-u$ and $w^0=U-u_0$. Then this difference $w$ satisfies the mild PDE
\begin{equation}\label{milddifference:PDE}
w(t)= e^{t\Delta}w^0 -\int_0^te^{(t-s)\Delta}\mathbb{P}\nabla\cdot (u\otimes w+w\otimes U)(s) \,ds.
\end{equation}
Moreover, our smallness assumptions on $u_0$ and $f$ and the usual fixed point Lemma~\ref{Can}
imply that $w$ can be obtained
as the limit in $L^\infty(\R^+,L^{3,\infty})$ of the approximating sequence $(w_k)$,
defined by
\begin{equation*}
 w_{k+1}=e^{t\Delta}w^0-\tilde B(u,w_{k})-\tilde B(w_k,U),
\end{equation*}
where the recursive relation starts with $w_0(x,t)=e^{t\Delta}w^0$.
Moreover, this sequence $(w_k)$ is bounded in $L^\infty(\R^+,L^{p,\infty})$.

By the semigroup property (recall that $F(x,t)$ denotes the kernel of $e^{t\Delta}\P\hbox{div}$):
$$\tilde B(u,v)(t)=e^{t\Delta/2}\tilde B(u,v)(t/2)+\int_{t/2}^t F(t-s)*(u\otimes v)(s)\,ds.$$
We deduce,
\begin{multline*}
w_{k+1}(t) =e^{t\Delta}w^0 -e^{t\Delta/2}\tilde B(u,w_k)(t/2)-e^{t\Delta/2}\tilde B(w_k,U)(t/2)\\
-\int_{t/2}^tF(t-s)*(u\otimes w_k)(s)\,ds
-\int_{t/2}^t F(t-s)*(w_k\otimes U)(s)\,ds.
\end{multline*}

Now let $r=\min(3,p)$ and denote
\begin{equation*}
M= \max \{ \|w^0\|_{L^{r,\infty}}, \|U\|_{L^{r,\infty}}, \sup_{s>0}\|u(s)\|_{L^{r,\infty}}\}.
\end{equation*}

By Lemma~\ref{heatforce:lemma} and using that the sequence $w_k$ is bounded in $L^\infty(\R^+,L^{3,\infty})$,
$$ \|\tilde B(u,w_k)(t/2)\|_{L^{r,\infty}} +\|\tilde B(w_k,U)(t/2)\|_{L^{r,\infty}} 
      \le C_rM. $$
A heat kernel estimate now implies, for all $q>r$ and for some constant $C'_r>0$ independent of~$q$,
\begin{multline*}
 \|w_{k+1}(t)\|_{L^{q,\infty}}
\le C'_{r}Mt^{-\frac{3}{2}(\frac{1}{r}-\frac{1}{q})} 
+\biggl\| \int_{t/2}^t F(t-s)*(u\otimes w_k)(s)\,ds\biggr\|_{L^{q,\infty}}\\
\qquad
+\biggl\| \int_{t/2}^t  F(t-s)*(w_k\otimes U)(s)\,ds\biggr\|_{L^{q,\infty}}.
\end{multline*}

From Lemma~\ref{heatforce:lemma} with H\"older's inequality
we have
\begin{multline*}
\biggl\| \int_{t/2}^t F(t-s)*(u\otimes w_k)(s)\,ds\biggr\|_{L^{q,\infty}} + \biggl\| \int_{t/2}^t  F(t-s)*(w_k\otimes U)(s)\,ds\biggr\|_{L^{q,\infty}}\\
\le C''_q\,\varepsilon_3\, \sup_{s\in [t/2,t]} \|w_k(s)\|_{L^{q,\infty}}.
\end{multline*}

Let
$$ W_{k}(t)\equiv\sup_{\tau\in[t,\infty)} \|w_{k}(\tau)\|_{L^{q,\infty}}. $$
then
\begin{equation*}
 W_{k+1}(t)\le C'_{r}M t^{-\frac{3}{2}(\frac{1}{r}-\frac{1}{q})} + C''_{q}\,\varepsilon_3\,W_k(t/2).
\end{equation*}

Iterating this inequality implies
\begin{align*}
W_k(t)
&\le C'_{r}M \sum_{n=0}^{k-1}
\Bigl(C''_{q}\,\varepsilon_3\,2^{\frac{3}{2}(\frac{1}{r}-\frac{1}{q})}\Bigr)^{n}
   t^{-\frac{3}{2}(\frac{1}{r}-\frac{1}{q})} +\bigl(C''_{q}\,\varepsilon_3\bigr)^k W_0(t/2^k)\\
&\le 2C'_{r}Mt^{-\frac{3}{2}(\frac{1}{r}-\frac{1}{q})}
    +C(r,q)\bigl(C''_{q}\,\varepsilon_3 2^{\frac{3}{2}(\frac{1}{r}-\frac{1}{q})}\bigr)^k t^{-\frac{3}{2}(\frac{1}{r}-\frac{1}{q})},
\end{align*}
provided
\begin{equation*}
C''_{q}\,\varepsilon_3\,2^{\frac{3}{2}(\frac{1}{r}-\frac{1}{q})}<\frac12 .
\end{equation*}
A slightly more stringent smallness condition and independent on $r>\frac{3}{2}$ is, e.g.,
\begin{equation}
\label{qdi8}
\varepsilon_3<\varepsilon(q):= \frac1{4C''_q}
\end{equation}

Now assuming~\eqref{qdi8} and letting $k\to\infty$ we get,
\begin{equation}
 \|w(t)\|_{L^{q,\infty}}\le 2C'_{r}M t^{-\frac{3}{2}(\frac{1}{r}-\frac{1}{q})}, \qquad \text{for } r=\min(3,p) \text{ and } q>r.
\end{equation}
Writing the above estimate for $q-\eta$ and $q+\eta$, for some $\eta>0$ small enough and interpolating
the $L^q$-space between $L^{q-\eta,\infty}$ and $L^{q+\eta,\infty}$ shows 
that the above estimate remains valid  with $\|w(t)\|_{L^q}$ on the left-hand side.
This establishes the stability result~\eqref{strong-q3p}.

\medskip

We now prove Claim (ii).
The weak convergence $u(t)\rightharpoonup U$ in $L^{p,\infty}$ for $\frac{3}{2}<p\le 3$
is obvious since the solution $u(t)$ is bounded in $L^{p,\infty}$
 and goes to $U$ in the sense of distributions  (even in $L^q$, $q>3$, as implied by the previous part of the proof).

On the other hand,  the proof of the necessary and sufficient condition for the strong convergence result in the $L^{3,\infty}$-norm is given
in \cite[Theorem 2.2]{CanK05} and in \cite[Corollary 4.1]{CanK05}, hence we will skip it.
The necessary and sufficient condition for the strong convergence result in the $L^{p,\infty}$-norm, with $\frac{3}{2}<p<3$
is a direct consequence of Claim (iii) which we now prove.

\medskip

Let now $\frac{3}{2}<p<3$ and  $\frac{3p}{6-p}\leq q\leq p$.
Given these restrictions, there exists some $\qb\in[p,4]$ such that the following relations hold true:
\begin{equation*}
\frac1q-\frac1p\leq\frac1{\qb}\leq \frac1p,\qquad \frac1{\qb}<\min\bigl(1-\frac1p,\frac13+\frac1q-\frac1p\bigr).  
\end{equation*}

We go back to the equation for $w$ given in \eqref{milddifference:PDE}.
We estimate $\|w(t)-w_0(t)\|_{L^q}$ using Propositions~\ref{holder} and~\ref{young},
the bound $\|w(t)\|_{L^{\qb,\infty}}\le Ct^{-\frac{3}{2}(\frac{1}{p}-\frac{1}{\qb})}$ (consequence of~\eqref{strong-q3p} with $q=\qb$)
and the the fact that $U\in L^{p,\infty}$ and $u\in L^\infty((0,\infty),L^{p,\infty})$.
We get
\begin{align*}
 \|w(t)-w_0(t)\|_{L^q}
&\le C\int_0^t (t-s)^{-\frac12+\frac32(\frac1q-\frac1p-\frac1{\qb})}\|(u\otimes w+w\otimes U)(s)\|_{L^{\frac{\qb p}{p+\qb},\infty}}\,ds\\
&\le C\int_0^t (t-s)^{-\frac12+\frac32(\frac1q-\frac1p-\frac1{\qb})}\|w(s)\|_{L^{\qb,\infty}}\bigl(\|u(s)\|_{L^{p,\infty}}+ \|U\|_{L^{p,\infty}}\bigr)\,ds\\
&\le C\int_0^t (t-s)^{-\frac12+\frac32(\frac1q-\frac1p-\frac1{\qb})}s^{-\frac{3}{2}(\frac{1}{p}-\frac{1}{\qb})}\,ds\\
&\le Ct^{\frac12+\frac3{2q}-\frac3p}.
\end{align*}
The theorem is now completely proved.

We finally show our stability result for large solutions. In lay words the next theorem shows that global solutions~$u$ of 
the non-stationary Navier--Stokes equations with a small constant forcing term~$f$ tend to forget their initial data.
Indeed, such solutions will converge to the small stationary solution~$U$ driven by the same forcing~$f$ even if they are initially large.

\begin{theorem}\label{largesol}
There exists an absolute constant $\eppatru>0$ with the following property. Let $u\in L^\infty_{loc}([0,\infty);L^{3,\infty})\cap L^4_{loc}([0,\infty);L^4)$ be a global solution of the evolutionary Navier-Stokes equations with a constant in time forcing $f$ such that $\Delta^{-1}\P f\in L^{3,\infty}\cap L^4$ and 
\begin{equation}\label{smallcondl3}
A(u_0,f)\equiv\limsup_{R\to0}R\mes\{|u_0|>R\}^{\frac13}+\|\Delta^{-1}\P f\|_{L^{3,\infty}}<\eppatru.  
\end{equation}
Let $U\in L^{3,\infty}\cap L^4$ be the unique stationary solution constructed in Theorem \ref{theoLp}.
Then we have that
\begin{itemize}
\item $\limsup\limits_{t\to\infty}\|u(t)\|_{L^{3,\infty}}\leq 22A(u_0,f);$
\item $u(t)\rightharpoonup U$ weakly in $L^{3,\infty}$ as $t\to\infty$;
\item $u(t)\to U$ in $L^{3,\infty}$ as $t\to\infty$ if and only if $e^{t\triangle}(u_0-U)\to 0$ strongly in $L^{3,\infty}$ as $t\to\infty$.
\end{itemize}
\end{theorem}
\begin{proof}
The idea of the proof is the same as in \cite{GIP} where it was proved that any global solution of the Navier--Stokes equations without
external force goes to 0 in the Besov spaces 
$\dot B^{-1+\frac3p}_{p,q}$ when the time becomes large (see also \cite{ADT} for the case of $VMO^{-1}$). It consists in decomposing the initial velocity in a small part plus a square integrable part. The small part remains small by the small data theory and the square-integrable part will become small at some point by using some energy estimates.

Here we use Lemma \ref{decomp} to decompose $u_0=v_0+w_0$ where $v_0\in L^2\cap L^{3,\infty}$ and $\|w_0\|_{L^{3,\infty}}<2A(u_0,f)$. Assuming that $3\eppatru<\eptrei$ where $\eptrei $ is the constant from Theorem \ref{forcedNS:theorem}, we can apply that theorem to construct a global solution $w$ of the Navier-Stokes equations with forcing term $f$, initial velocity $w_0$ and such that
\begin{equation*}
\|w(t)\|_{L^{3,\infty}} \leq 8A(u_0,f) \qquad \text{for all }t\geq0.
\end{equation*}
Moreover, according to relation \eqref{strong-q3p} the solution $w$ satisfies the following decay estimate $\sup \limits_{t>0}t^{\frac18}\|w(t)-U\|_{L^4}<\infty$. Since $U\in L^4$ we infer that  $w\in L^4_{loc}([0,\infty);L^4)$.

The difference $v=u-w$ verifies the following PDE:
\begin{equation}\label{eqv}
\partial_tv-\Delta v+u\cdot \nabla v+v\cdot\nabla w+\nabla p'=0  
\end{equation}
whose integral form reads
\begin{equation}\label{integv}
v(t)=e^{t\Delta}v_0-\int_0^t e^{(t-s)\Delta}\P\nabla\cdot(u\otimes v+v\otimes w)(s)\,ds.  
\end{equation}

We show first that $v\in C^0([0,\infty);L^2)$. The first term on the RHS above clearly belongs to this space. We show that so does the second term. The kernel $F(t)$ of the operator $e^{t\Delta}\P\operatorname{div}$ is of the form $F(x,t)=t^{-2}F(\frac x{\sqrt t},1)$ with $F(\cdot,1)\in L^1\cap L^\infty\subset L^{p,q}$ for all $1<p<\infty$ and $1\leq q\leq\infty$. In particular, $\|F(t)\|_{L^{\frac65,2}}\leq t^{-\frac34}$ so that $F\in L_{loc}^1([0,\infty);L^{\frac65,2})$. By the H\"older inequality we also have that $u\otimes v+v\otimes w\in L_{loc}^\infty([0,\infty);L^{\frac32,\infty})$. Since the last term in \eqref{integv} is the space-time convolution of $F$ with   $u\otimes v+v\otimes w$, we infer that it belongs to  $C^0([0,\infty);L^2)$.

For $0<\delta<1$, let $J_\delta$ be a smoothing operator that multiplies in the frequency space by a cut-off function bounded by 1 which is
a smoothed out version of the characteristic function of the annulus $\{\delta<|\xi|<\frac1\delta\}$. We also introduce an approximation of the identity
$\varphi_\eta$ in time. 

Given the additional regularity found for $v$ above, we remark that we can multiply the equation of $v$ expressed in \eqref{eqv} by $\varphi_\eta\ast\varphi_\eta \ast J_\delta^2v$ and integrate in space and time from $t_0$ to $t$, with $t_0>0$, to obtain that
\begin{multline}\label{estvinit}
\|\varphi_\eta\ast J_\delta v(t)\|_{L^2}^2+2\int_{t_0}^t \|\nabla\varphi_\eta\ast J_\delta v(s)\|_{L^2}^2\,ds
=\|\varphi_\eta\ast J_\delta v(t_0)\|_{L^2}^2 \\+2\int_0^t\int_{\R^3}u\cdot\nabla (\varphi_\eta\ast \varphi_\eta\ast J^2_\delta v)\cdot v
+2\int_0^t\int_{\R^3}v\cdot\nabla (\varphi_\eta\ast \varphi_\eta\ast J^2_\delta v)\cdot w.  
\end{multline}
We let now $\eta\to0$. Given the time continuity of $v$ with values in $L^2$, we have that $\varphi_\eta\ast J_\delta v(t)\to J_\delta v(t)$ and $\varphi_\eta\ast J_\delta v(t_0)\to J_\delta v(t_0)$ in $L^2$ as $\eta\to0$. The other terms in  \eqref{estvinit} pass easily to the limit $\eta\to0$. Therefore, taking first the limit $\eta\to0$ in \eqref{estvinit}, and second $t_0\to0$ and using again that $v\in C^0[0,\infty);L^2)$ we get that 
\begin{multline}\label{estv}
\|J_\delta v(t)\|_{L^2}^2+2\int_0^t \|\nabla J_\delta v(s)\|_{L^2}^2\,ds
=\|J_\delta v_0\|_{L^2}^2 +2\int_0^t\int_{\R^3}u\cdot\nabla J^2_\delta v\cdot v\\
+2\int_0^t\int_{\R^3}v\cdot\nabla J^2_\delta v\cdot w
\end{multline}

We bound the last two terms on the RHS as follows
\begin{multline*}
2\int_0^t\int_{\R^3}u\cdot\nabla J^2_\delta v\cdot v
+2\int_0^t\int_{\R^3}v\cdot\nabla J^2_\delta v\cdot w  
\leq 2\int_0^t \|\nabla J^2_\delta v\|_{L^2}\|v\|_{L^4}(\|u\|_{L^4}+\|w\|_{L^4})\\
\leq \frac12  \int_0^t\|\nabla J_\delta v\|^2_{L^2} +  \int_0^t \|v\|^2_{L^4}(\|u\|^2_{L^4}+\|w\|^2_{L^4}).
\end{multline*}
Plugging this in \eqref{estv} yields
\begin{equation*}
\|J_\delta v(t)\|_{L^2}^2+\int_0^t \|\nabla J_\delta v(s)\|_{L^2}^2\,ds
\leq \|J_\delta v_0\|_{L^2}^2 +\int_0^t \|v\|^2_{L^4}(\|u\|^2_{L^4}+\|w\|^2_{L^4}).
\end{equation*}
Since $u,v,w\in L^4_{loc}([0,\infty);L^4)$, the RHS above is uniformly bounded with respect to $\delta$. Letting $\delta\to0$ implies thanks to the Beppo-Levi theorem that $ \int_0^t \|\nabla v(s)\|_{L^2}^2\,ds<\infty$, that is $v\in L^2_{loc}([0,\infty);H^1)$. 

We go back to \eqref{estv} and estimate
\begin{multline}\label{v1}
2\int_0^t\int_{\R^3}u\cdot\nabla J^2_\delta v\cdot v
=2\int_0^t\int_{\R^3}u\cdot\nabla J^2_\delta v\cdot (1-J_\delta^2)v\\
\leq 2\|u\|_{L^4(0,t;L^4)} \|\nabla v\|_{L^2(0,t;L^2)} \|(1-J_\delta^2)v\|_{L^4(0,t;L^4)} \to0 \quad\text{as }\delta\to0. 
\end{multline}

We observe now that $\dot H^1(\R^3)\hookrightarrow L^{6,2}(\R^3)$. This imbedding follows from the Young inequality for Lorentz spaces after noticing that $(-\Delta)^{-\frac12}$ is a convolution operator with a function bounded by $\frac C{|x|^2}$ which therefore belongs to $L^{\frac32,\infty}$. We use this fact together with the H\"older inequality to bound the last term in \eqref{estv} as follows
\begin{equation}\label{v2}
2\int_0^t\int_{\R^3}v\cdot\nabla J^2_\delta v\cdot w
\leq C \int_0^t\|v\|_{L^{6,2}} \|\nabla J^2_\delta v\|_{L^2}\|w\|_{L^{3,\infty}} 
\leq CA(u_0,f) \int_0^t\|\nabla v\|^2_{L^2} 
\end{equation}
Using \eqref{v1} and \eqref{v2} in \eqref{estv}, letting $\delta\to0$ and using the Beppo-Levi theorem we infer that 
\begin{equation*}
\|v(t)\|_{L^2}^2+2\int_0^t \|\nabla v(s)\|_{L^2}^2\,ds
\leq\|v_0\|_{L^2}^2 +CA(u_0,f) \int_0^t\|\nabla v\|^2_{L^2}
\end{equation*}
If we further assume that $C\eppatru\leq1$, then $CA(u_0,f)\leq1$ so the relation above implies that $v\in L^\infty(\R_+;L^2)\cap L^2(\R_+;\dot H^1)$. By interpolation and from the imbedding $\dot H^{\frac12}\subset L^{3,\infty}$ we infer that $v\in L^4(\R_+;\dot H^{\frac12})\subset L^4(\R_+; L^{3,\infty}) $. So there exists a time $T=T(\eppatru)$ such that $\|v(T)\|_{L^{3,\infty}}< A(u_0,f)$. Since we also have that $\|w(T)\|_{L^{3,\infty}}< 8 A(u_0,f)$ we infer that $\|u(T)\|_{L^{3,\infty}}< 9 A(u_0,f)$. Assuming further that $10\eppatru<\eptrei$, Theorem \ref{forcedNS:theorem} allows to construct a small solution starting from time $T$, a solution whose $L^{3,\infty}$ norm will be bounded by $22 A(u_0,f)$. We will prove below a uniqueness result stating that $u$ must be equal to this small solution starting from time $T$. Once this is proved, the first part of the theorem follows. Moreover, using again that our solution $u$ becomes small after the time $T$, the second and the third part of the theorem are consequences of Theorem \ref{forcedNS:theorem}. Except that the equivalent condition for $u(t)$ to converge strongly to $U$ in $L^{3,\infty}$ is that $e^{t\triangle}(u(T)-U)\to 0$ strongly in $L^{3,\infty}$ as $t\to\infty$. To finish the proof it therefore suffices to show that
\begin{equation*}
e^{t\triangle}(u(T)-U)\stackrel{t\to\infty}{\longrightarrow} 0 
\text{ in }L^{3,\infty}\quad \Longleftrightarrow \quad  e^{t\triangle}(u_0-U)\stackrel{t\to\infty}{\longrightarrow} 0 
\text{ in }L^{3,\infty}.
\end{equation*}
This is a consequence of the following sequence of equivalence relations:
\begin{align*}
e^{t\triangle}(u(T)-U)\stackrel{t\to\infty}{\longrightarrow} 0 
\text{ in }L^{3,\infty}\ 
&\Longleftrightarrow \  e^{t\triangle}(w(T)-U)\stackrel{t\to\infty}{\longrightarrow} 0 
\text{ in }L^{3,\infty}\\
&\Longleftrightarrow \  w(t)\stackrel{t\to\infty}{\longrightarrow} U 
\text{ in }L^{3,\infty}\\
&\Longleftrightarrow \  e^{t\triangle}(w_0-U)\stackrel{t\to\infty}{\longrightarrow} 0 
\text{ in }L^{3,\infty}\\
&\Longleftrightarrow \  e^{t\triangle}(u_0-U)\stackrel{t\to\infty}{\longrightarrow} 0 
\text{ in }L^{3,\infty}.
\end{align*}
We used above that $v(T),v_0\in L^2$ and the decay estimates for the heat equation to deduce the first and fourth lines of the relation above, and Theorem \ref{forcedNS:theorem} twice for $w$, starting from time $t=0$ and from time $t=T$ to deduce the second and third lines. This completes the proof provided that we prove the announced uniqueness result.

Let $\ub$ be the small solution starting from time $T$ with initial velocity $u(T)$ constructed in Theorem \ref{forcedNS:theorem} and set $\vb=\ub-w$. As above we have that $\vb\in C^0([T;\infty);L^2)\cap L^2([T;\infty);\dot H^1)$. Then $v-\vb$ solves the following equation:
\begin{equation*}
\partial_t(v-\vb)-\Delta(v-\vb)+u\cdot\nabla(v-\vb)+(v-\vb)\cdot\nabla\ub=-\nabla p_1.  
\end{equation*}
As in the previous argument, one can prove that this relation can be multiplied by $v-\vb$ and integrated from $T$ to $t$ to get that, for all $t\geq T$, 
\begin{multline*}
  \|(v-\vb)(t)\|_{L^2}^2+2\int_{T}^t\|\nabla(v-\vb)\|_{L^2}^2
=\int_{T}^t\int(v-\vb)\cdot\nabla(v-\vb)\cdot\ub\\
\leq C\int_{T}^t \|\nabla(v-\vb)\|_{L^2}^2\|\ub\|_{L^{3,\infty}}
\leq CA(u_0,f)\int_{T}^t \|\nabla(v-\vb)\|_{L^2}^2
\leq \int_{T}^t \|\nabla(v-\vb)\|_{L^2}^2
\end{multline*}
provided that $A(u_0,f)$ is sufficiently small. We infer that $v(t)=\vb(t)$, that is $u(t)=\ub(t)$ for all $t\geq T$. This completes the proof of the theorem.
\end{proof}

\begin{remark}
We also have stability in $L^{p,\infty}$ for large solutions. More precisely, suppose that in addition to the hypothesis of Theorem \ref{largesol}  we assume that $u_0\in L^{p,\infty}$ with $p\in(\frac32,3)$. Then $u\in L^\infty(\R_+;L^{p,\infty})$,  $u(t)\rightharpoonup U$ weakly in $L^{p,\infty}$ and    $u(t)-e^{t\Delta}(u_0-U)\to U$ in $L^{p,\infty}$ as $t\to\infty$. This follows easily after applying Theorem \ref{forcedNS:theorem} starting from the time $T$ when  the solution becomes small. One only needs to show the following two facts:
\begin{itemize}
\item if $u_0\in L^{p,\infty}$ and $\Delta^{-1}\P f\in L^{p,\infty}$ then $u\in L^\infty(0,T;L^{p,\infty})$;
\item $e^{(t-T)\Delta}(u(T)-U)-u^{t\Delta}(u_0-U)\to 0$ strongly in $L^{p,\infty}$ as $t\to\infty$.
\end{itemize}
To prove the first assertion, we observe that, with the notation from the proof of Theorem \ref{forcedNS:theorem} (namely the notation used in relation \eqref{mildNS:PDE2}) one has that $\tilde u_0\in L^\infty(\R_+;L^{p,\infty})$. Moreover, by the H\"older inequality and using the standard decay estimates for the heat equation we can bound
\begin{multline}\label{Best}
\|\tilde B(u,u)(t)\|_{L^{p}}\leq
\int_0^t \|e^{(t-s)\Delta}\P\nabla\cdot(u\otimes u)(s)\|_{L^{p}}\,ds\\
\leq C\int_0^t(t-s)^{-\frac32+\frac3{2p}} \|u(s)\|^2_{L^{3,\infty}}
\leq Ct^{-\frac12+\frac3{2p}}\sup_{0<s<t}\|u(s)\|^2_{L^{3,\infty}}.
\end{multline}
We infer that $\tilde B(u,u)\in L^\infty(0,T;L^{p})\subset L^\infty(0,T;L^{p,\infty})$, so by \eqref{mildNS:PDE2} we also have that $u\in L^\infty(0,T;L^{p,\infty})$. To show the second assertion, we observe that it is sufficient to prove that
\begin{equation*}
u(T)-U-e^{t\Delta}(u_0-U)\in L^{q,\infty}  
\end{equation*}
for some $q<p$. But $u-U$ verifies the PDE
\begin{equation*}
\partial_t(u-U)-\Delta(u-U)+u\cdot\nabla u-U\cdot\nabla U=-\nabla p_2  
\end{equation*}
whose mild formulation implies that 
\begin{equation*}
u(T)-U-e^{t\Delta}(u_0-U)=-\int_0^Te^{(T-s)\Delta}\P\nabla\cdot(u\otimes u-U\otimes U)(s)\,ds
\end{equation*}
The same estimate as in \eqref{Best} shows now that the RHS belongs to $L^{q,\infty}$ for any $\frac32<q<3$, in particular for some $q<p$.

Moreover, if $u_0\in L^p$, then the previous argument shows that $u\in L^\infty(0,T;L^p)$. From Theorem \ref{forcedNS:theorem} applied starting from time $T$ we infer that  $u\in L^\infty(\R_+;L^p)$ and $u(t)\to U$ in $L^p$ as $t\to\infty$.  
\end{remark}

\begin{remark}
We observe that the condition imposed on the initial velocity by the hypothesis of Theorem \ref{largesol} does not imply that $u_0$ is close in $L^{3,\infty}$ to the smooth functions in $\mathcal{S}(\R^3)$.
 Indeed, that would require to have that the quantity $\limsup\limits_{R\to\infty}R\mes\{|u_0|>R\}^{\frac13}$ is small too. 
This condition is not necessary
in Theorem \ref{largesol}.
\end{remark}

\section{Acknowledgements}

Theorem~\ref{theoLp} is a development of an insightful remark made to the first and the last author
by an anonymous referee of their paper~\cite{BjorSch07}.
The authors gratefully acknowledge him.
The authors would like to thank also the referees of the present journal for their many useful suggestions.

\bibliographystyle{amsplain}

\end{document}